\documentclass{article}
\usepackage{amsmath,amssymb,amsthm,textcomp,color}
\usepackage[all]{xy}

\usepackage[utf8]{inputenc}
\usepackage[breaklinks]{hyperref}
\theoremstyle{plain}
\newtheorem{theorem}{Theorem}
\newtheorem{corollary}[theorem]{Corollary}
\newtheorem{proposition}[theorem]{Proposition}
\newtheorem{lemma}[theorem]{Lemma}
\newtheorem{example}{Example}
\newtheorem{remark}{Remark}

\theoremstyle{definition}
\newtheorem{definition}[theorem]{Definition}

\voffset = - 1.2cm
\numberwithin{equation}{section}
\numberwithin{theorem}{section}
\input xy
\xyoption{all}
\def\di{\bigstar}

\numberwithin{equation}{section}
\numberwithin{theorem}{section}

\begin{document}

\centerline{{\Large \bf Quasi-Lie families, schemes, invariants}}
\vskip 0.25cm
\centerline{{\Large \bf and their applications to Abel equations}}
\vskip 0.5cm
\centerline{J.F. Cari\~nena$^\dagger$ and J. de Lucas$^\ddagger$}
\vskip 0.25cm

\centerline{$^\dagger$ Departamento de  F\'{\i}sica Te\'orica and IUMA, Universidad de
Zaragoza,}
\vskip 0.1cm
\centerline{c. Pedro Cerbuna 12, 50009 Zaragoza, Spain.}
\centerline{$^\ddagger$ Department of Mathematical Methods in Physics, University of Warsaw}
\vskip 0.10cm
\centerline{ul. Pasteura 5, 02-093, Warszawa, Poland}
\vskip 0.25cm

\begin{abstract}

This work analyses types of group actions
on families of $t$-dependent vector fields of a particular class, the
hereby called {\it quasi-Lie families}. We devise methods to obtain the defined here {\it quasi-Lie invariants}, namely a kind of functions constant along
the orbits of the above-mentioned actions.
Our techniques lead to a deep
geometrical understanding of quasi-Lie schemes and quasi-Lie systems giving rise to several new results.
Our achievements are
illustrated by studying Abel and Riccati equations. We retrieve the Liouville invariant and study other new quasi-Lie invariants of
Abel equations. Several Abel equations with a superposition rule are described and we characterise Abel equations via quasi-Lie schemes.

\end{abstract}

\section{Introduction}
The theory of Lie systems investigates a class of non-autonomous systems of first-order ordinary differential
equations, the {\it Lie systems}, admitting a {\it superposition rule},
i.e. a function permitting us to describe the general solution
of such a system in terms of
a generic family of particular solutions and some constants  \cite{CGM07,Dissertationes,LS,PW}. Lie systems
have been attracting
attention due to their geometric properties and applications \cite{Br95,Dissertationes}.

Despite 
providing useful geometric methods to study Lie systems and
possessing many applications, the theory of Lie systems presents an important drawback:
just a few, but important,
systems of differential equations can be studied through this theory \cite{Dissertationes,In72}. To employ the theory of Lie systems for analysing more general systems, the notion of {\it quasi-Lie scheme} 
 emerged \cite{CGL08}. This concept led to
the description of new and known properties of non-linear oscillators
\cite{CGL08}, second-order Gambier equations
\cite{CGL13}, Emden--Fowler equations \cite{CLL08Emd}, second-order Riccati equations \cite{CLSecRicc}, dissipative Milne--Pinney
equations \cite{Dissertationes},   Abel equations \cite{CLRAbel}, etc.

Roughly speaking, quasi-Lie schemes describe the transformation properties of certain
collections of non-autonomous systems of first-order ordinary differential equations of a particular kind, the so-called {\it quasi-Lie families},
under the action of a related group of $t$-dependent changes of variables
\cite[Theorem 1]{CGL08}. As a simple application, this enables us to map Abel equations onto Abel equations providing methods to integrate them \cite{CLRAbel}.

Many procedures developed to study the integrability of a non-autonomous system of differential equations of a certain sort, e.g. Abel equations, rely on mapping it onto
an autonomous and easily integrable one of the same type by transformations
generally provided in an  {\it ad hoc}\, way. Meanwhile,
the theory of quasi-Lie schemes may naturally
provide answers to these questions \cite{CLRAbel}, e.g. it justifies geometrically the group of transformations employed to study Abel equations \cite{CLRAbel}. 
Since such answers depend only on the
algebraic/geometric structure of the system of differential equations under study, the so found properties can easily  be generalised to other systems of differential equations with a similar structure.

As a first new result, 
we succeed in characterising groups of $t$-dependent
transformations mapping elements of a quasi-Lie family onto other (possibly the same)  members of the quasi-Lie family: the so-called {\it structure preserving groups} of the quasi-Lie family. 
The elements of such a group of transformations, the Lie almost symmetries, can be considered as a generalisation of the well-known Lie symmetries. Our new theory recovers as particular 
cases the (extended) quasi-Lie schemes groups for quasi-Lie schemes described in \cite{CGL08,CLL08Emd}. We also succeeded in describing the whole family of $t$-dependent transformations 
mapping Abel equations onto Abel equations adding details not discussed in \cite{CLRAbel}.

Next, we characterise Abel equations as elements of the quasi-Lie family associated with a quasi-Lie scheme. This shows that many of the properties of
Abel equations are essentially algebraic and geometric and are due to the structure of a quasi-Lie scheme. This motivates to define new structures related to quasi-Lie schemes: 
quasi-Lie schemes morphisms, the representation associated with a quasi-Lie scheme, etc. On the one hand, these notions enlarge the theoretical content of the theory of quasi-Lie schemes. On the other hand, 
we generalise the so-called hierarchies of Lie systems given in  \cite{BM10} and we extrapolate the results  on  integrability of Abel equations to other equations, e.g. Lotka--Volterra systems, and vice versa.

Another important question is whether we can connect two elements of a quasi-Lie family by means of a Lie almost symmetry. To provide necessary conditions for the existence of such a Lie almost symmetry, we define the {\it quasi-Lie invariants} of a quasi-Lie family, i.e. functions mapping each element of the quasi-Lie family onto a $t$-dependent function that is the same for all elements connected by a Lie almost symmetry.

We describe and analyse different types of quasi-Lie invariants and  we characterise them in geometrical and algebraics terms. As an application, we retrieve the known Liouville invariant of Abel equations and
 we prove the existence of other new quasi-Lie invariants.

The study of Abel equations is currently a very active field of research that aims, for example, to determinate their integrability conditions, exact solutions, and
general properties \cite{Al07,Ch40,Kamke,MCH01,Pa05,SP03,Stre,SR82,W12}. For instance, Lie symmetries of Abel equations were studied in \cite{I14}. Their solutions via transformations $y=u(x)z(x)+v(x)$ were 
 investigated in \cite{SMSL} and through the referred to as Julia's condition in \cite{Bou14}; equivalence and integrable cases are presented in \cite{CR03}; etc. 
 Meanwhile, some solutions in closed form appeared in \cite{NSST} and periodic solutions have been studied in \cite{AGG06}. Abel equations are also relevant for other 
 mathematical problems like the centre problem for planar polynomial systems of differential equations \cite{Al13,AB14,GGS13,GGS14,P14}.
Moreover, Abel equations are broadly used in Physics, where they
describe the properties of a number of interesting physical systems
\cite{CK86,CLP01,AMM97,GIM96,HM04,MML08,Pa05,SP03,ZT09}. For instance, Abel equations play a role  in the analysis of cosmological models \cite{HM03,HM04,YYY}
 and the theory of membranes in $M$-theory \cite{CLP01,ZT09}. Hence, the analysis
of Abel equations is very relevant  so as to determine the properties of the physical systems they describe.

The structure of our paper goes as follows. Section 2 contains the fundamental notions and results of the theory of Lie and quasi-Lie
schemes that are used throughout this work. Sections 3 and 4 introduce structure preserving groups and Lie almost symmetries for quasi-Lie families. 
As a consequence, some results about the transformation  
properties of Riccati and Abel equations are retrieved.
Sections 5 and 6 suggest new algebraic and geometric constructions for the study of quasi-Lie schemes and the systems of differential equations investigated by means of them.
 In Section 7 we study Abel equations through quasi-Lie schemes. Section 8 is devoted to the study of integrable Abel equations. We introduce and investigate quasi-Lie invariants in Section 9. 
 Subsequently, Section 10 is aimed to study quasi-Lie invariants of order two. Our results in previous section are employed to analyse Abel equations in Section 11. 
 Section 12 is devoted to studying quasi-Lie invariants of order one and zero and their use in the study of Abel equations. In Section 13 we summarise our results and sketch our future work. 
 
\section{Fundamentals notions on Lie systems and quasi-Lie schemes}

Let us describe the fundamental results of the
theories of Lie systems and quasi-Lie schemes used in this work (see \cite{CGL08,CGM07,Dissertationes} for details). To stress the main points of our presentation, we focus on systems of differential equations on vector spaces and assume all mathematical objects to be global and smooth. A rigorous theory on manifolds can easily be developed from our results.

Given the projection $\pi:(t,x)\in \mathbb{R}\times\mathbb{R}^n\mapsto
x\in\mathbb{R}^n$ and the tangent bundle projection $\tau:{\rm
T}\mathbb{R}^n\rightarrow\mathbb{R}^n$, a {\it $t$-dependent vector field} on
$\mathbb{R}^n$ is a mapping $X:\mathbb{R}\times\mathbb{R}^n\rightarrow{\rm
T}\mathbb{R}^n$ such that $\tau\circ X=\pi$. Hence, every
$t$-dependent vector field $X$ is equivalent to a family $\{X_t\}_{t\in\mathbb{R}}$ of
vector fields $X_t:x\in \mathbb{R}^n\mapsto X(t,x)\in{\rm T}\mathbb{R}^n$.

An {\it integral curve} of $X$ is a standard
integral curve
$\gamma:\mathbb{R}\rightarrow\mathbb{R}\times \mathbb{R}^n$ of its {\it suspension}, i.e. the vector field $\overline{ X}=\partial/\partial t+X(t,x)$ on $\mathbb{R}\times\mathbb{R}^n$.

Given a non-autonomous system of differential equations
\begin{equation}\label{Sys}
\frac{dx^i}{dt}=X^i(t,x),\qquad i=1,\ldots,n,
\end{equation}
we can define
a unique $t$-dependent vector field on $\mathbb{R}^n$, namely
\begin{equation}\label{TdVF}
X(t,x)=\sum_{i=1}^nX^i(t,x)\frac{\partial}{\partial x^i},
\end{equation}
whose integral curves of the form $(t,x(t))$ are the particular solutions to (\ref{Sys}). Conversely, given a generic $t$-dependent vector field (\ref{TdVF}) there exists a unique system
of first-order differential equations, namely (\ref{Sys}), whose particular solutions are the integral curves of $X$ of the form $\gamma:t\in\mathbb{R}\mapsto (t,x(t))\in\mathbb{R}\times\mathbb{R}^n$. This is
the referred to as {\it associated system} with $X$. This justifies to write $X$ for both a
$t$-dependent vector field and its associated system.

Let us denote by $\mathfrak{X}(\mathbb{R}^n)$ and $\mathfrak{X}_t(\mathbb{R}^n)$ the spaces of vector fields and $t$-dependent vector fields on
$\mathbb{R}^n$, respectively. Each $X\in\mathfrak{X}_t(\mathbb{R}^n)$ gives rise to a {\it generalised flow} $g^X$, i.e. a mapping
$g^X:(t,x)\in\mathbb{R}\times \mathbb{R}^n \mapsto  g^X_t(x)\equiv
g^X(t,x)\in \mathbb{R}^n$, with $g^X_0={\rm Id}_{\mathbb{R}^n}$, such that the
curve $\gamma^X_{x_0}(t)=g^X_t(x_0)$ is the particular solution of the system associated with $X$ with initial condition $\gamma^X_{x_0}(0)=x_0$.
We call {\it autonomisation} of $g^X$ the diffeomorphism
$\overline{g}^X:(t,x)\in\mathbb{R}\times\mathbb{R}^n\mapsto
(t,g^X_t(x))\in\mathbb{R}\times\mathbb{R}^n$. 

Every generalised flow $h^X$ (related to a generic element  $X\in \mathfrak{X}_t(\mathbb{R}^n)$) acts on $\mathfrak{X}_t(\mathbb{R}^n)$ transforming each $Y\in\mathfrak{X}_t(\mathbb{R}^n)$ into a new one $h^X_\di Y\in\mathfrak{X}_t(\mathbb{R}^n)$ with
generalised flow $(g^{h_\di Y})_t=h^X_t\circ g^Y_t$ for all $t\in\mathbb{R}$. We can
write (see \cite[Theorem 3]{CGL08})
\begin{align*}
\overline{h^X_\di Y}=\bar{h}^X_*\overline{Y}\,,
\end{align*}where $\bar{h}^X_*$ is the standard action of the diffeomorphism $\bar{h}^X$ on
vector fields of $\mathfrak{X}(\mathbb{R}\times \mathbb{R}^n)$.

A {\it Lie system} is a non-autonomous system of first-order ordinary
differential equations admitting a {\it superposition rule}, i.e. a map
$\Phi:{\mathbb{R}}^{nm}\times\mathbb{R}^n\to {\mathbb{R}}^n$, $x=\Phi(x_{(1)},
\ldots,x_{(m)};k_1,\ldots,k_n)$, such that the general solution, $x(t)$,  of the system
can be written as 
\begin{align*}x(t)=\Phi(x_{(1)}(t),
\ldots,x_{(m)}(t);k_1,\ldots,k_n),\end{align*} 
where $\{x_{(a)}(t)\mid
a=1,\ldots,m\}$ is any `generic' family of particular solutions and
$k=(k_1,\ldots,k_n)$ is a set of $n$ constants \cite{CGLTdep,Dissertationes}.

Superposition rules do not depend explicitly on $t$, which in physical motivated examples usually refers to the time. Hence,
superposition rules only remain well defined under certain, very general, changes of variables, e.g. those of the form $\phi: (t,x)\in \mathbb{R}\times \mathbb{R}^n\mapsto
(\bar t\equiv \phi_1(t),\bar x\equiv \phi_2(x))\in \mathbb{R}\times \mathbb{R}^n$ for diffeomorphisms $\phi_1$ and $\phi_2$. Indeed, a $t$-dependent change of variables $\bar x=\varphi(t,x(t))$ maps the initial Lie system onto a new one
which does not need to admit a superposition rule \cite{CGL12}. Moreover, the $t$-independence of
superposition rules is the key point to find their many properties, e.g. the Lie--Scheffers theorem (see \cite{CGM07}).

Lie characterised systems (\ref{Sys}) admitting a superposition rule \cite{LS}. His result, the Lie--Scheffers theorem \cite{CGM07}, states that a system $X$ admits
a superposition rule if and only if
\begin{equation}
X(t,x)=\sum_{\alpha =1}^r b_\alpha(t)\, X_{(\alpha )}(x),
\label{Lievf}
\end{equation}
for certain $t$-dependent functions $b_1(t),\ldots,b_r(t)$ and a set  $\{X_{(\alpha)}\mid \alpha=1,\ldots,r\}$  of vector fields generating an
$r$-dimensional real Lie algebra: a {\it
Vessiot--Guldberg Lie algebra} of $X$. In other words, $X$ is a Lie system if and only if $\{X_t\}_{t\in\mathbb{R}}$ is contained in a
finite-dimensional real Lie algebra of vector fields.
The Lie--Scheffers theorem suggests us to define the following structure, which also occurs  naturally when defining quasi-Lie schemes \cite{CGL08}.

\begin{definition}
We call {\it quasi-Lie family} the $C^\infty(\mathbb{R})$-module $\mathcal{V}$ of $t$-dependent vector
fields $X\in\mathfrak{X}_t(\mathbb{R}^n)$ taking values in a finite-dimensional vector space $V$ of vector fields on $\mathbb{R}^n$.
\end{definition}


\begin{definition} A {\em quasi-Lie scheme} $S(W,V)$ consists of two
finite-dimensional vector spaces of vector fields
$W,V\subset\mathfrak{X}(\mathbb{R}^n)$ such that 
 $W$ is a Lie algebra of vector fields
contained in  $V$ and $W$ normalises $V$, that is, $[W,W]\subset W\subset V$
and  $[W,V]\subset V$. We call {\it group of } $S(W,V)$ the group $(\mathcal{G}_W,\star)$ of generalised flows of $t$-dependent vector fields of $\mathcal{W}$ with the product 
\begin{align*}
(g^X\star h^Y)_t=g^X_t \circ h^Y_t,\qquad \forall t\in\mathbb{R},\qquad \forall g^X,h^Y\in \mathcal{G}_W.
\end{align*}
\end{definition}

Given a quasi-Lie scheme
$S(W,V)$, the quasi-Lie family  $\mathcal{V}$ is `stable' under the action of
the elements of $\mathcal{G}_W$, i.e. $g^X_\di Y\in
\mathcal{V}$, for every $Y\in \mathcal{V}$ and $g^X\in \mathcal{G}_W$ (see  \cite[Theorem 1]{CGL08}).

Consider a Lie group action  $\Phi:G\times \mathbb{R}^n\rightarrow \mathbb{R}^n$.  The Lie algebra of fundamental vector fields of the action $\Phi$ will be denoted $W_\Phi$. The diffeomorphisms
$\Phi_g:x\in\mathbb{R}^n\mapsto\Phi(g,x)\in\mathbb{R}^n, \forall g\in G$,
 act on 
$\mathfrak{X}_t(\mathbb{R}^n)$ transforming every
$X\in\mathfrak{X}_t(\mathbb{R}^n)$ into $X'\in \mathfrak{X}_t(\mathbb{R}^n)$ given by
$X'_t=\Phi_{g*}X_t$. To keep notation simple, we use $g$ instead of $\Phi_g$ when the meaning is clear.

Additionally, every curve $g:t\in\mathbb{R}\mapsto g_t\in G$ transforms 
$X\in\mathfrak{X}_t(\mathbb{R}^n)$ into a new one $X'\in\mathfrak{X}_t(\mathbb{R}^n)$ via the 
diffeomorphism
$\bar g:(t,x)\in\mathbb{R}\times\mathbb{R}^n\mapsto
(t,\Phi(g_{t},x))\in\mathbb{R}\times \mathbb{R}^n$, i.e.
 $\bar g_{*}\bar X=\bar X'$ for
a unique $X'\in\mathfrak{X}_t(\mathbb{R}^n)$. In short, the curve $g$ in $G$
transforms every $X\in\mathfrak{X}_t(\mathbb{R})$ into a
new one $X'=g_\di X$. Moreover, it  was shown in \cite{CGL08} that
\begin{equation}\label{Descom}
(g_\di X)_t=\frac{\partial g_t}{\partial t}\circ g_t^{-1}+g_{t*}X_t,\qquad \forall t\in\mathbb{R}.
\end{equation}

The finite-dimensional Lie algebra $W$ of $S(W,V)$ can be integrated to give rise to a local Lie group action 
$\Phi:G_W\times \mathbb{R}^n\rightarrow  \mathbb{R}^n$ with $T_eG\simeq W$.  The Lie group action becomes global when the vector fields of $W_\Phi$ are complete  (cf. \cite{Dissertationes}).

\section{$t$-independent structure preserving groups}
The Riccati equation is a Lie system on the real line related to a Vessiot--Guldberg Lie algebra $V_{\rm Ricc}$
isomorphic to $\mathfrak{sl}(2,\mathbb{R})$ \cite{PW}. Once that a one-point compactification of the real line is carried out, then becoming a circle $S^1$, 
the $V_{\rm Ricc}$ can be integrated to give rise to a Lie group action $\Phi:PSL(2,\mathbb{R})\times S^1\rightarrow S^1$ whose maps
$\Phi_g$, with $g\in PSL(2,\mathbb{R})$, transform every Riccati equation $X\in \mathcal{V}_{\rm Ricc}$ onto a new Riccati equation $\Phi_{g\star} X\in \mathcal{V}_{\rm Ricc}$ (cf. \cite{CarRam}). 
In this section we aim to generalise and to study this result to quasi-Lie families.

Each diffeomorphism $\Phi_g$ can
be understood as
a generalisation of the symmetry notion for Riccati equations that, instead of mapping particular solutions of a Riccati equation onto particular solutions of the same Riccati equation,
transforms solutions of a Riccati equation into solutions of another Riccati equation. This suggests us	 the following definition.

\begin{definition}
A Lie group action
$\Phi:G\times\mathbb{R}^n\rightarrow\mathbb{R}^n$ is a {\it $t$-independent structure preserving group of transformations} of the quasi-Lie family
$\mathcal{V}$ if for every $g\in G$ and $X\in \mathcal{V}$, we have $g_\star X\in \mathcal{V}$. In this case, we call each $\Phi_g$ a {\it Lie almost symmetry}
of $X\in \mathcal{V}$.
\end{definition}
To shorten the terminology, we frequently  say  structure preserving group instead of structure preserving group of transformations. Now, it is easy to prove the following result.
\begin{lemma} Let $X,Y\in \mathfrak{X}(\mathbb{R}^n)$ and let
$h^Y:\mathbb{R}\times\mathbb{R}^n\rightarrow\mathbb{R}^n$ be the flow
of $Y$, then
\begin{align*}
\frac{\partial}{\partial s}h^{Y}_{-s*}X=h^Y_{-s*}[Y,X],\qquad \forall s\in\mathbb{R}.
\end{align*}
Moreover,
\begin{equation}\label{PrimDecom}
(h^Y_{s*}
X)_x=X_x-s[Y,X]_x+\frac{s^2}{2!}[Y,[Y,X]]_x-\frac{s^3}{3!}[Y,[Y,[Y,X]]]_x+\ldots
\end{equation}
for every $x\in \mathbb{R}^n$ where the right-hand side of the
above expression converges.
\end{lemma}

A Lie group action
$\Phi:G\times\mathbb{R}^n\rightarrow\mathbb{R}^n$, with $G$ being connected, is a $t$-independent structure preserving  group of a quasi-Lie family
$\mathcal{V}$ if and only if $[W_\Phi,V]\subset V$. Indeed,  if $\Phi$ is a $t$-independent structure preserving group for $\mathcal{V}$, then every element of $V$ is stable under the flow of every $Y\in W_\Phi$ and $[Y,V]\subset V$. Conversely, if $[W_\Phi,V]\subset V$, then $g_\star X\in  \mathcal{V}$
 for every element $g\in G$  and each $X\in\mathcal{V}$.

Let us use the above property to enlighten results
given by the theory of Lie systems. Consider the set of Riccati  differential equations
\begin{equation}\label{Ricc}
\frac{dx}{dt}=a_0(t)+a_1(t)x+a_2(t)x^2,
\end{equation}
with $a_0(t),a_1(t),a_2(t)$ being arbitrary $t$-dependent functions.  Each such equation  is described by a $t$-dependent vector field $X\in \mathcal{V}_{{\rm
Ricc}}\subset \mathfrak{X}_t(\mathbb{R})$, with
\begin{equation}\label{VGRicc}
V_{{\rm Ricc}}=\left\langle\frac{\partial}{\partial x},x\frac{\partial}{\partial
x},x^2\frac{\partial}{\partial x}\right\rangle.
\end{equation}
It can be proved that the (local) Lie group action
$\Phi:U\subset PSL(2,\mathbb{R})\times \mathbb{R}\rightarrow\mathbb{R}$, where $PSL(2,\mathbb{R})=SL(2,\mathbb{R})/\{I,-I\}$,  for $I$  the $2\times 2$ identity  matrix $I$,
\begin{align*}
\Phi(g,x)=
{\frac{\alpha x+\beta}{\gamma x+\delta}},
\qquad {\rm where} \,\, g=\left[\left(\begin{array}{cc}\alpha&\beta\\
\gamma&\delta\end{array}\right)\right]\in PSL(2,\mathbb{R}),
\end{align*}
i.e. $g$ is the equivalence class of a matrix of $SL(2,\mathbb{R})$, the set $U=\{(g,x)\in PSL(2,\mathbb{R})\times \mathbb{R}\,|\, x\neq -\delta/\gamma\}$,  is such that for every
$g\in PSL(2,\mathbb{R})$ and $X\in \mathcal{V}_{{\rm Ricc}}$, we have $g_\star X\in
\mathcal{V}_{{\rm Ricc}}$ \cite{CarRam}. In other words, $\Phi$ is a $t$-independent structure preserving group of $\mathcal{V}_{{\rm Ricc}}$ and
 $[W_\Phi,V_{{\rm Ricc}}]\subset
V_{{\rm Ricc}}$. Indeed, it is easy to verify that $W_\Phi=V_{{\rm Ricc}}$
(cf. \cite{Br95}). Now, a natural question arises: is there
any other $t$-independent structure preserving group $\Phi_2:\widetilde{G}\times \mathbb{R}\rightarrow\mathbb{R}$, with $PSL(2,\mathbb{R})\subset \widetilde{G}$ being connected and $\Phi_2$ 
being effective, up to a  zero measure subset of $\widetilde G$, 
for $\mathcal{V}_{\rm Ricc}$?

To answer the previous question, let us determine every Lie algebra $W$ that satisfy that $[W,V_{{\rm
Ricc}}]\subset V_{{\rm Ricc}}$. Given any $X_1=f(x)\partial/\partial x
\in W$, the condition $[W,V_{\rm Ricc}]\subset V_{\rm Ricc}$ entails
\begin{align*}
\left[X_1,\frac{\partial}{\partial x}\right]=-f'(x)\frac{\partial}{\partial x}\in V_{{\rm Ricc}},
\end{align*}
where $'$ stands for the derivative in terms of $x$. So, $f(x)=\sum_{\alpha=0}^3c_\alpha x^\alpha$ for certain $c_0,c_1,c_2,c_3\in \mathbb{R}$. From this and using that $[X_1,x^2\partial/\partial x]\in V_{\rm Ricc}$, we obtain
\begin{align*}
\left[X_1,x^2\frac{\partial}{\partial x}\right]
=(2x\,f(x)-x^2\,f'(x))\frac{\partial}{\partial x}=(2c_0x+c_1\,x^2-c_3\,x^4)\frac{\partial}{\partial x}\in V_{\rm Ricc}.
\end{align*}
Hence, $c_3=0$. Consequently, the largest Lie
algebra $W$ with $[W,V_{\rm Ricc}]\subset V_{\rm Ricc}]$, i.e.
preserving $\mathcal{V}_{{\rm Ricc}}$, is $W=V_{{\rm Ricc}}$. This fixes the possible
Lie algebras of fundamental vector fields of $t$-independent structure preserving groups for Riccati equations and the form close to the identity of the
connected Lie group of the corresponding action.

Likewise, we  can  analyse the first-order {\it Abel equations of
the first kind and degree $q$}, namely
\begin{equation}\label{Abel}
\frac{dx}{dt}=f_0(t)+f_1(t)x+\ldots+f_q(t)x^q,\qquad q\geq 3,
\end{equation}
where $f_0(t), \ldots, f_{q}(t)$, are arbitrary $t$-dependent functions \cite{Abel,Al07}. 

\begin{remark}Our definition of Abel equations does not require $f_q$ to be different from the zero function as it appears in other works. 
Nevertheless, we can still prove similar results to those given in our work for these latter Abel equations by introducing minor modifications.
\end{remark}

Abel equations
 are related to  $t$-dependent vector fields of the
quasi-Lie family $\mathcal{V}_{{\rm Abel}}$ with
\begin{equation}\label{VectAbel}
V_{{\rm Abel}}=\left\langle Y_0,Y_1,\ldots,Y_q\right\rangle,\qquad Y_k=x^k\frac{\partial}{\partial x},\qquad k=0,\ldots,q.
\end{equation}

We aim to determine  a $t$-independent structure preserving group, i.e.  the `largest' finite-dimensional Lie algebra of fundamental vector fields $W$ satisfying $[W,V_{{\rm Abel}}]\subset V_{{\rm Abel}}$. 
Given an arbitrary $X_1\in W$, let us say $X_1=f(x)\partial/\partial x$, the condition
\begin{align*}
\left[X_1,Y_0\right]=-f'(x)\frac{\partial}{\partial x}\in V_{{\rm
Abel}},
\end{align*}
implies that $f(x)$  must be a polynomial with a  degree $p\leq q+1$.
Moreover, $X_1$ must satisfy
\begin{align*}
\left[X_1,Y_k\right]=(f(x)kx^{k-1}-f'(x)x^k)\frac{\partial}{\partial x}\equiv g_k(x)\frac{\partial}{\partial x}\in
V_{{\rm Abel}},\qquad k=1,\ldots,q.
\end{align*}
Observe that ${\rm deg}\, g_k(x)=p+k-1$ for $k\neq p$ and ${\rm deg}\, g_p(x)<2p-1$. Since $[W,V_{\rm Abel}]\subset V_{\rm Abel}$, then ${\rm deg}\, g_k(x)\leq q$ for every $k$ and the only possibilities are $p=0,1$. Hence, the largest Lie algebra $W$ satisfying $[W,V_{\rm Abel}]\subset V_{\rm Abel}$ is
\begin{align*}
W_{{\rm Abel}}=\left\langle\frac{\partial}{\partial x},x\frac{\partial}{\partial x}\right\rangle.
\end{align*}
Integrating this Lie algebra, we obtain  a  structure preserving group for  (\ref{VectAbel}):
\begin{equation}\Phi_{\rm Abel}:(\beta,\alpha;x)\in (\mathbb{R}\rtimes\mathbb{R}^+)\times \mathbb{R}\mapsto \bar x\in\mathbb{R},\qquad \bar x\equiv \alpha\, x+\beta,\label{transf}
\end{equation}
where $\mathbb{R}\rtimes\mathbb{R}^+$ denotes the semi-direct product Lie group with composition law
\begin{align}\label{ActAbe}
(\beta,\alpha)\star (\beta',\alpha')\equiv (\alpha\beta'+\beta,\alpha \alpha'),\qquad (\beta,\alpha),(\beta',\alpha')\in \mathbb{R}\rtimes\mathbb{R}^+.
\end{align}

The above Lie group action of $\mathbb{R}\rtimes\mathbb{R}^+$ on $\mathbb{R}$ can be considered as a particular case of a new $t$-independent structure preserving group 
of $\mathcal{V}_{\rm Abel}$ based on $\mathbb{R}\rtimes\mathbb{R}_*$, where $\mathbb{R}_*\equiv\mathbb{R}-\{0\}$. To do so, we 
must find Lie almost symmetries $\varphi:x\in \mathbb{R}\mapsto \bar x\in \mathbb{R}$ of Abel equations which are not related to the previous structure preserving group. 
Each Lie almost symmetry for Abel equations obeys that $\varphi_*$ maps elements of $V_{\rm Abel}$ onto elements of $V_{\rm Abel}$. Moreover, $\varphi_*$ 
is a  $C^\infty(\mathbb{R})$-linear Lie algebra automorphism of $\mathfrak{X}(\mathbb{R})$. So, the value $\varphi_*Y$ for any $Y\in V_{\rm Abel}$ can be
 obtained from the value $\varphi_*Y_0$ with $Y_0=\partial/\partial x$. In view of the previous facts and since the linear map ${\rm ad}_{Y_0}:Y\in V_{\rm Abel}\mapsto [Y_0,Y]\in V_{\rm Abel}$ is nilpotent, 
 then ${\rm ad}_{\varphi_*Y_0}$ must also be a nilpotent endomorphism when acting on $V_{\rm Abel}$. Thus, the only possibility is $\varphi_*Y_0=k_0Y_0$, with $k_0\in \mathbb{R}_*$,
  namely $\varphi:x\in\mathbb{R}\mapsto k_0x+k_1\in\mathbb{R}$ for certain $(k_1,k_0)\in\mathbb{R}\rtimes\mathbb{R}_*$.  When $k_0<0$, the previous transformation is not one
   of the Lie almost symmetries of $\Phi_{\rm Abel}$. Hence, we can enlarge $\Phi_{\rm Abel}$ to a Lie group action of $G_{W_{\rm Abel}}\equiv\mathbb{R}\rtimes \mathbb{R}_*$ with 
   composition law (\ref{ActAbe}).
This completely characterises the set of affine transformations employed to study all Abel equations \cite{ES04} and fulfills the analysis given in \cite{CLRAbel}.

\section{Structure preserving groups}
The $t$-independent structure preserving groups are a natural particular case of a more general object, the hereafter defined structure preserving groups, that enable us to describe geometrically the transformations appearing in the
theory of integrability of non-autonomous differential equations and 
 Lie systems. More specifically,  the set  $\mathcal{G}$ of curves in  a Lie group $G$  is a group when endowed with  the composition law:
\begin{align*}
(g\star h)_t=g_t\cdot h_t,\qquad g,h \in\mathcal{G},
\end{align*}
where `$\cdot$' denotes the composition law of $G$. The
 neutral element of $(\mathcal{G},\star)$ is the constant curve $t\mapsto e$, with $e$ being the neutral element of $G$.

We showed in \cite{CGL08} that every curve in $G$ acts on the space of $t$-dependent vector fields on $\mathbb{R}^n$ and  given two curves $g$ and $h$ in $G$, we have
${g}_\di({h}_\di X)=(g\star h)_\di X$. From this,  every Lie group action
$\Phi:G\times\mathbb{R}^n\rightarrow\mathbb{R}^n$ gives rise to a Lie group action
$\Theta:\mathcal{G}\times\mathfrak{X}_t(\mathbb{R}^n)\rightarrow
\mathfrak{X}_t(\mathbb{R}^n)$.

\begin{definition} A Lie group action
$\Phi:G\times\mathbb{R}^n\rightarrow\mathbb{R}^n$ is a {\it
structure preserving group of (transformations of) a quasi-Lie family
 $\mathcal{V}$} if $g_\di X\in \mathcal{V}$ for every curve $g$ in $G$ and
each $X\in \mathcal{V}$.  We call $g_\di$ a {\it Lie almost symmetry} of $\mathcal{V}$.
\end{definition}
As before, we will simply  say structure preserving group to simplify the notation, and in full similarity with the $t$-independent case, we can prove that 
an action
$\Phi:G\times\mathbb{R}^n\rightarrow\mathbb{R}^n$, with $G$ being a connected group,
is a structure preserving group of a quasi-Lie family $\mathcal{V}$ if and
only if the pair $(W_\Phi, V)$ is a quasi-Lie
scheme.

Consequently,  every $S(W,V)$
 defines a unique action, the hereby called {\it action of the quasi-Lie scheme},
 $\Phi:G_W\times\mathbb{R}^n\rightarrow\mathbb{R}^n$, with $G_W$ being a simply connected Lie
group with $\dim\, G_W=\dim\, W$, whose fundamental vector fields are those in $W$. The
elements of the group $\mathcal{G}$ of curves in $G_W$ give rise to a group of
$t$-dependent transformations of $\mathbb{R}^n$ of the form
$g:(t,x)\in\mathbb{R}\times\mathbb{R}^n\mapsto \Phi(g_t,x)\in\mathbb{R}^n $
mapping elements of $\mathcal{V}$ onto elements of $\mathcal{V}$. 

The above construction can be extended to a Lie group action where $G_W$ is just connected. In this case, the group of curves in $G_W$ retrieves the
 so-called {\it extended group of the scheme}, ${\rm Ext}(W)$,
defined in \cite{CLL08Emd}. Moreover, if we assume $g_0=e$, the previous
transformations give rise to the generalised flows of $t$-dependent vector
fields of $W$, i.e. the group $\mathcal{G}_W$, introduced for the first time in
\cite{CGL08}.

Let us illustrate the results of the previous section through the
analysis of  Abel equation of the first kind for $q=3$. Nevertheless, results can straightforwardly be extended
to $q>3$.

As $S(W_{{\rm Abel}},V_{{\rm Abel}})$ is a quasi-Lie scheme, then  $\Phi_{\rm Abel}$ given by (\ref{transf})
 is a structure preserving group of transformations for Abel equations. That is, the $t$-dependent changes of variables induced by
curves in $\mathbb{R}\rtimes\mathbb{R}_*$, namely $x(t)=\alpha(t)\bar x(t)+\beta(t)$ with
$\alpha(t)\neq 0$ for every $t\in\mathbb{R}$, transform each 
$X\in \mathcal{V}_{{\rm Abel}}$ into a new $t$-dependent vector field of $ \mathcal{V}_{{\rm Abel}}$. Indeed,
given a curve $g^{-1}:t\in\mathbb{R}\mapsto (\beta(t),\alpha(t))\subset \mathbb{R}\rtimes\mathbb{R}_+$ and
$X=\sum_{\alpha=0}^3f_\alpha (t)Y_\alpha$, it can be proved that
$
g_\di X=\sum_{\alpha=0}^3\bar{f}_\alpha(t)Y_\alpha\in \mathcal{V}_{{\rm Abel}},
$
where 
\begin{equation}\label{NewCoff}
\begin{aligned}
\bar{f}_3(t)&=f_3(t)\alpha^2(t),\\
\bar{f}_2(t)&=\alpha(t)(f_2(t)+3f_3(t)\beta(t)),\\
\bar{f}_1(t)&=3f_3(t)\beta^2(t)+2f_2(t)\beta(t)+f_1(t)-\frac{1}{\alpha(t)
}\frac{d\alpha}{dt}(t),\\
\bar{f}_0(t)&=\frac{1}{\alpha(t)}
\left(f_3(t)\beta^3(t)+f_2(t)\beta^2(t)+f_1(t)\beta(t)+f_0(t)-\frac{d\beta}{dt}(t)\right).\\
\end{aligned}
\end{equation}
Therefore, the integral curves of $g_\di X$ of the form $(t,\bar x(t))$ are particular solutions of
\begin{equation}\label{finalAbel}
\frac{d\bar x}{dt}=\bar{f}_0(t)+\bar{f}_1(t)\bar x+\bar{f}_2(t)\bar x^2+\bar{f}_3(t)\bar{x}^3.
\end{equation}

As  there is a group
of $t$-dependent transformations  related to each scheme, we can
generalise the concept of quasi-Lie system given in \cite{CGL08} as follows:

\begin{definition} A system $X\in \mathcal{V}$  is a {\it quasi-Lie system} for $S(W,V)$ and a  Lie group $G_W$ if there is a curve $g$ in $G_{W}$ such that
$g_\di X\in \mathcal{V}_0\subset \mathcal{V}$, for a Lie algebra $V_0$ of vector fields.
\end{definition}

If $X$ is a quasi-Lie system with respect to $S(W,V)$, then the geometric properties of $g_\di X$, e.g. its superposition rules,
can be studied through the theory of Lie systems. Posteriorly, the $t$-dependent
change of variables induced by $g$ can be used to recover the properties of $X$. Superposition rules for $g_\di X$ need not give rise to superposition rules
for $X$ because the change of variables induced by $g$ may be $t$-dependent.

\section{Quasi-Lie families, representations and characterisation of Abel equations}

Given the quasi-Lie family $\mathcal{V}$ related to $S(W,V)$, we can make use of $S(W,V)$ to define a new structure, the representation associated with the scheme, which is a key to study $\mathcal{V}$.
As an application, we characterise Abel equations of order $q\geq 3$ as elements of the quasi-Lie family $\mathcal{V}$ related to a quasi-Lie scheme $S(W,V)$ on $\mathbb{R}$ with $\dim W=2$ and $\dim V=q+1$.
\begin{definition}
We call {\it representation associated with $S(W,V)$} the Lie algebra morphism
\begin{align*}
\begin{array}{rccc}
\rho_{W,V}:& W &\longrightarrow &{\rm End}(V)\\
                  &X_1&    \mapsto     &{\rm ad}_{X_1}
\end{array}
\end{align*}
where ${\rm ad}_{X_1}: X_2\in V\mapsto [X_1,X_2]\in V$ and ${\rm End}(V)$, i.e.  the space of linear endomorphisms on $V$, is endowed with the Lie bracket given by the commutator of endomorphisms on $V$.
\end{definition}

\begin{remark}
Observe that $\rho_{W,V}$ is well defined: since $[W,V]\subset V$, then ${\rm ad}_{X_1}\in {\rm End}(V)$ for every $X_1\in W$.
This leads to a Lie algebra representation of $W$ on the space ${\rm End}(V)$. Although $\rho_{W,V}$ is an algebraic concept,
 the elements of $W$ and $V$ are vector fields, which implies that this representation is also geometrical. 
 It is also worth noting that the general theory of integration of Lie algebra homomorphisms proves that  $\rho_{W,V}$ determines the action of the structure preserving  group 
$\Phi:G_W\times \mathbb{R}^n\rightarrow \mathbb{R}^n$ on $\mathcal{V}$.
\end{remark}


Let us use above ideas to characterise Abel equations through  quasi-Lie schemes.

\begin{theorem}\label{ClaAbel}
A first-order ordinary differential equation  on the real line in normal form is an Abel equation of degree $q\geq 3$  (up to a local $t$-independent change of variables) if and only if it can be considered as an element of the quasi-Lie family $\mathcal{V}$ of a quasi-Lie scheme $S(W,V)$ on $\mathbb{R}$ with $\dim W=2$ and $\dim V=q+1$.
\end{theorem}
\begin{proof}

The direct part is a trivial consequence of our previous results. Let us prove the converse.
Since we assume $\dim W=2$, then $W$ must be either Abelian or isomorphic to $\mathfrak{Affin}(\mathbb{R})$, namely the Lie algebra of vector fields corresponding to the Lie group action of affine transformations on the real line. Every Lie algebra of vector fields on the real line is locally isomorphic to a Lie subalgebra of $\mathfrak{sl}(2,\mathbb{R})$. Hence,  $W$ cannot be Abelian and it admits a basis $X_1,X_2$ obeying that $[X_1,X_2]=X_1$ \cite{Hum}. We can easily find  a certain coordinate system on $\mathbb{R}$ driving a certain basis $\{X_1,X_2\}$ of $W$ with $[X_1,X_2]=X_1$ into the form $X_1=\partial/\partial x$ and $X_2=x\partial/\partial x$.

  We can extend ${\rm ad}_{X_1}$ to a unique $\mathbb{C}$-linear mapping on the whole complexified vector space $V_\mathbb{C}\equiv\mathbb{C}\otimes V$. Since $V_\mathbb{C}$ is a vector space over $\mathbb{C}$, then ${\rm ad}_{X_1}$ has at least an eigenvector, namely, a vector field $Z_1\in V_\mathbb{C}\backslash\{0\}$ satisfying that
$$
[X_1,Z_1]=\lambda_0 Z_1\Rightarrow Z_1=\lambda_1 e^{\lambda_0 x}\frac{\partial}{\partial x},
$$
 for a certain $\lambda_0\in\mathbb{C}$ and $\lambda_1\in\mathbb{C}\backslash\{0\}$. We fix $\lambda_1=1$ for simplicity. As $X_1$ is a real vector field and $V_\mathbb{C}$ is closed under complex conjugation, we have that $\bar Z_1=e^{\bar \lambda_0 x}\partial/\partial x\in V_\mathbb{C}$ satisfies  that
$$
[X_1,\bar Z_1]=\bar{\lambda}_0 \bar Z_1.
$$
As a consequence,  $Z_+\equiv (Z_1+\bar Z_1)/2=e^{{\rm Re}(\lambda_0)x}\cos({\rm Im}(\lambda_0) x)\partial/\partial x$ and also $Z_-\equiv (Z_1-\bar Z_1)/(2i)=e^{{\rm Re}(\lambda_0)x}\sin({\rm Im}(\lambda_0) x)\partial/\partial x$ belong to $V$.  All the Lie brackets of ${\rm ad}_{X_2}^nZ_+$ and ${\rm ad}_{X_2}^nZ_-$ belong to $V$ if and only
$$
{\rm ad}_{X_2}^n\left(e^{{\rm Re}(\lambda_0)x}e^{i{\rm Im}(\lambda_0)x}\frac{\partial}{\partial x}\right)\in V_\mathbb{C},\qquad \forall n\in\mathbb{N}.
$$
For $\lambda_0\neq 0$ the above vector fields generate an infinite-dimensional family of vector fields linearly independent over $\mathbb{C}$ belonging to $V_\mathbb{C}$. As $V$ is finite-dimensional, this is impossible and $\lambda_0=0$. Since ${\rm ad}_{X_1}$ is triangularizable on $V_\mathbb{C}$ and its unique eigenvalue is zero, it must be a nilpotent endomorphism on $V_\mathbb{C}$ and, as $X_1$ is a real vector field, on $V$. Hence, $V$ takes the form (\ref{VectAbel}) 
for a certain $q\in \mathbb{N}$ with $q\geq 3$. This shows that every quasi-Lie scheme on the real line with $\dim W=2$ and $\dim\,V\geq 4$ must be of the above form $S(W,V)$ up to a change of variables.

Since our $t$-dependent vector field $X$ takes values in $V$, it must take the form of an Abel equation of degree $p$ (up to a $t$-independent change of variables).
\end{proof}

\section{Algebraic morphisms of quasi-Lie schemes and Abel equations}

Theorem \ref{ClaAbel} suggests that the properties of Abel equations must be encoded in the structure $S(W_{\rm Abel},V_{\rm Abel})$. Likewise, systems related to quasi-Lie schemes sharing the same structure will presumably have similar properties. This suggests us to define the following.

\begin{definition} An {\it algebraic quasi-Lie schemes morphism} between two quasi-Lie schemes $S(W_1,V_1)$ and $S(W_2,V_2)$  is a linear mapping $\phi:V_1\rightarrow V_2$ such that:
\begin{itemize}
 \item ${\rm ad}_{\phi(X_1)}\phi(Y_1)=\phi({\rm ad}_{X_1}Y_1)$ for every $X_1\in W_1$ and $Y_1\in V_1$,
 \item $\phi(W_1)\subset W_2$.
\end{itemize}
If $\phi(W_1)=W_2$ and $\phi(V_1)=V_2$, we say that $\phi$ is an {\it algebraic quasi-Lie schemes epimorphism}. If $\phi$ is a linear monomorphism, we say that $\phi$ is {\it an algebraic quasi-Lie schemes monomorphism}. We say that $\phi$ is an {\it algebraic quasi-Lie isomorphism} if it is an algebraic monomorphism and epimorphism of quasi-Lie schemes.
\end{definition}

The above definitions do not take into account the geometric nature of quasi-Lie schemes. That is why we call them algebraic. Despite that, they allow us to describe relevant  transformation properties of quasi-Lie schemes. In view of previous definitions, we can prove quite easily the following result.

\begin{proposition} If $\phi:V_1\rightarrow V_2$ is a quasi-Lie schemes morphism between $S(W_1,V_1)$ and $S(W_2,V_2)$, then $S(\phi(W_1),\phi(V_1))$ and $S(\ker \phi\cap W,\ker \phi)$ are quasi-Lie schemes.
\end{proposition}
\begin{minipage}[t]{3.5cm}
\begin{center}
\xymatrix{W_1\ar[d]^{\exp}\ar[r]^{\phi}&W_2\ar[d]^{\exp}\\
G_{W_1}\ar[r]^{\Phi}&G_{W_2}}
\end{center}
\end{minipage}
\begin{minipage}[t]{12cm}
Observe that a quasi-Lie schemes epimorphism $\phi:V_1\rightarrow V_2$ between $S(W_1,V_1)$ and $S(W_2,V_2)$ induces a Lie algebra epimorphism $\phi:W_1\rightarrow W_2$. Since $W_1$ and $W_2$ can be understood as the Lie algebras of $G_{W_1}$ and $G_{W_2}$, which are now assumed to be simply connected, we can define a unique Lie group epimorphism $\Phi$ making commutative the diagram on the left. Using the above, we can prove the following theorem.\\
\end{minipage}

\begin{theorem}\label{Sq} Let $\phi:V_1\rightarrow V_2$ be an algebraic epimorphism of quasi-Lie schemes $S(W_1,V_1)$  on $N$ and $S_2(W_2,V_2)$ on $N'$.
If $X\in \mathcal{V}_1$ is such that $g_{1\di} X=0$ for a certain curve $g_1$ in $G_{W_1}$, then
$\Phi (g_{1})_\di \phi (X)=0$, where $[\phi(X)]_t=\phi(X_t)$ for all $t\in\mathbb{R}$.
\end{theorem}
If $X'=g_{\di}X=0$ and $x(t)$ is a particular solution to $X$, then $g(x(t))$ is a particular solution to $X'$, namely $g(x(t))=x_0\in N$. Hence, $x(t)=g^{-1}(x_0)$ is a particular solution to $X$ and $g^{-1}$ is, essentially, the generalised flow associated with $X$. According to Theorem \ref{Sq}, we have that $\Phi(g)_\di\phi(X)=0$ and $\Phi(g)^{-1}$ is essentially the generalised flow of $\phi(X)$.  Additionally, we can also prove that if $X$ is a quasi-Lie system for $S(W_1,V_1)$, then $\phi(X)$ is a quasi-Lie system for $S(W_2,V_2)$. The above results retrieve as particular cases the findings about hierarchies of Lie systems in \cite{BM10}.

Let us show an example of our previous results. Consider the system
\begin{align}\label{Lot}
\frac{d\xi}{dt}=a_0(t)Z_0(\xi)+a_1(t)Z_1(\xi)+a_2(t)Z_2(\xi)+ a_3(t)Z_3(\xi),
\end{align}
where $\xi=(x,y)\in\mathbb{R}^2$ and
\begin{equation*}
Z_0=\frac{\partial}{\partial x},\quad Z_1=x\frac{\partial}{\partial x}+y\frac{\partial}{\partial y},\quad Z_2=(x^2-y^2)\frac{\partial}{\partial x}+2xy\frac{\partial}{\partial y},\quad Z_3=(x^3-3y^2x)\frac{\partial}{\partial x}+3x^2y\frac{\partial}{\partial y}.
\end{equation*}
This system of differential equations describes several cases of polynomial differential equations
appearing in the XVI Hilbert's problem \cite{Li03}. Moreover, the study of planar polynomial equations has also being attracting some attention \cite{Ma94,PW13}. By defining $S(W_2,V_2)$ with $W_2=\langle Z_0,Z_1\rangle$ and $V_2=\langle Z_0,\ldots,Z_3\rangle$, we can establish an algebraic isomorphism of quasi-Lie schemes $\phi:V_{\rm Abel}\rightarrow V_2$. In view of this and Theorem \ref{Sq}, every quasi-Lie system for $S(W_{\rm Abel},V_{\rm Abel})$, namely a certain Abel equation, is equivalent to a quasi-Lie system for $S(W_2,V_2)$.
\section{Quasi-Lie systems and Abel equations}
It is a well-known fact that two generic Abel equations cannot be connected through a $t$-dependent linear
transformation of the form $x(t)=\alpha(t)\bar x(t)+\beta(t)$, with $\alpha(t)\neq 0$ \cite{ES04}. This leads to the
theory of integrable classes introduced in the classical theory of integrability
of Abel equations. Roughly speaking, this theory analyses and uses the existence of a
transformation  of the above-mentioned form relating an initial Abel equation to an easier/integrable one to integrate the initial Abel equation \cite{ES04}.

Let us determine  Abel equations that are quasi-Lie systems with respect to
$S(W_{\rm Abel},V_{{\rm Abel}})$. More specifically, we search for  $t$-dependent vector
fields $X\in \mathcal{V}_{\rm Abel}$ admitting a curve $g$ in $G_{W_{\rm Abel}}$ such
that $g_\di X\in \mathcal{V}^{(2)}_\mu$, where $V^{(2)}_\mu$ is a Lie algebra spanned by the vector fields
\begin{align}\label{2mu}
Z_1=-2\mu^3Y_0+3\mu Y_2+Y_3,\qquad Z_2=\mu Y_0+Y_1.
\end{align}
The relevance of the Lie algebra $V^{(2)}_\mu$ is that every Lie algebra of vector fields in $V$ is either one-dimensional or it is of the form $V^{(2)}_\mu$ (cf. \cite{CLRAbel}). If $g_\di X\in \mathcal{V}^{(2)}_{\mu}$, then there exist
functions $\lambda_1(t)$ and $\lambda_2(t)$ such that $g_\di
X=\lambda_1(t)Z_1+\lambda_2(t)Z_2$. In order to determine $g$ in an easy way, it becomes
useful to write $g=\bar{g}^{-1}$ in terms of a curve $\bar g(t)=(\beta(t),\alpha(t))$, i.e. $x(t)=\alpha(t)\bar x(t)+\beta(t)$. In this manner, we obtain
\begin{align}
\lambda_1(t)\!&=\!f_3(t)\alpha^2(t),\label{eqv1}\\
3\mu\lambda_1(t)\!&=\!\alpha(t)(3f_3(t)\beta(t)+f_2(t)),\label{eqv2}\\
\lambda_2(t)\!&=\!3f_3(t)\beta^2(t)+2f_2(t)\beta(t)+f_1(t)-\frac{\dot\alpha(t)}{
\alpha(t)},\label{eqv3}\\
\alpha(t)(\mu\lambda_2(t)-2\mu^3\lambda_1(t))\!&=\!
f_3(t)\beta^3(t)\!+\!f_2(t)\beta^2(t)\!+\!f_1(t)\beta(t)\!+\!f_0(t)\!-\!\dot\beta(t)\label{eqv4}.
\end{align}
We hereafter assume for simplicity that $f_3(t)>0$ for every $t\in\mathbb{R}$. A similar procedure can be obtained for $f_3(t)<0$ for each $t\in\mathbb{R}$. The first equation of the above system implies that the sign of $\lambda_1(t)$ is the same as the sign of $f_3(t)$ at every $t\in\mathbb{R}$. Hence,
\begin{equation}\label{alpha}
\alpha(t)=\sqrt{\lambda_1(t)/f_3(t)},\qquad \forall t\in\mathbb{R},
\end{equation}
where, for the sake of simplicity, we choose $\alpha(t)>0$. Since $\alpha(t)\neq 0$ and $f_3(t)$, then $\lambda_1(t)\neq 0$ for every $t\in\mathbb{R}$. Therefore,
\begin{align*}
\frac{\dot \alpha}{\alpha}=\frac 12\left(\frac{\dot\lambda_1}{\lambda_1}-\frac{\dot f_3}{f_3}\right).
\end{align*}
Using these results in (\ref{eqv2})-(\ref{eqv4}), we obtain
\begin{equation}\label{sys3}
\begin{aligned}
3\mu\sqrt{\lambda_1}&=\frac{3f_3\beta+f_2}{\sqrt{f_3}},\\
\lambda_2&=3f_3\beta^2+2f_2\beta+f_1-\frac
12\left(\frac{\dot\lambda_1}{\lambda_1}-\frac{\dot
f_3}{f_3}\right),\\
\mu\lambda_2-2\mu^3\lambda_1&=\sqrt{\frac{f_3}{\lambda_1}}
\left(f_3\beta^3+f_2\beta^2+f_1\beta+f_0-\dot\beta\right
),
\end{aligned}
\end{equation}
where we recall that $\mu\in \mathbb{R}$. Assuming that $\mu\neq 0$, the first expression of the preceding system allows us to obtain $\lambda_1$
in terms of $\beta$, $\mu$ and the $t$-dependent coefficients of the initial
Abel equation. Moreover, above expressions and assumptions allow us to ensure that $f_3(f_2+3f_3\beta)\neq 0$. It follows that
\begin{equation}\label{l2}
\begin{gathered}
\lambda_2=f_1+2 f_2\beta+3f_3\beta^2+\frac{f_2\dot f_3-f_3\dot f_2-3f_3^2\dot \beta}{f_3(f_2+3f_3\beta)},\\
\lambda_2-\frac
2{9f_3}(3f_3\beta+f_2)^2=\frac{
3f_3\left(f_3\beta^3+f_2\beta^2+f_1\beta+f_0-\dot\beta\right)}{3f_3\beta+f_2}.\\
\end{gathered}
\end{equation}

Substituting $\lambda_2$ from the first expression into the second, we obtain
that
\begin{equation}\label{ConAbel}
9 f_3 \left(3 f_0 f_3+\dot f_2\right)-9 f_2 \left(f_1 f_3+\dot f_3\right)+2 f_2^3= 0.
\end{equation}
If this integrability condition is satisfied, we can obtain $\lambda_2$ by substituting any function $\beta\neq -f_2/(3f_3)$ in the first equation within (\ref{l2}). From $\beta$ and $\lambda_2$, the function $\lambda_1$ can be derived from the first equation in (\ref{sys3}). Then, $\alpha$ can be determined by means of (\ref{alpha}).

If we assume $\mu=0$, the first equation within (\ref{sys3}) shows that $\beta=-f_2/(3f_3)$. This and the third relation in (\ref{sys3}) imply that $\beta$ is a particular solution of the initial Abel equation. By using this, it can be proved that the coefficients of the Abel equation satisfy (\ref{ConAbel}). Moreover, the value of $\lambda_1(t)$ can be chosen an arbitrary positive function. Next,  $\alpha(t)$ can be obtained from (\ref{alpha}) and $\lambda_2(t)$ from (\ref{sys3}).

Summarising, we have proved the following remarkable theorem.
\begin{theorem}\label{th1}
An Abel equation of first-order and first-type with $f_3(t)>0$  is a quasi-Lie system with respect to the scheme $S(W_{\rm Abel},V_{\rm Abel})$ and a Lie algebra $V_\mu^{(2)}$ spanning by the vector fields (\ref{2mu}) if and only if
\begin{equation}\label{CA}
9 f_3 \left(3 f_0 f_3+\dot f_2\right)-9 f_2 \left(f_1 f_3+\dot f_3\right)+2 f_2^3= 0.
\end{equation}
\end{theorem}
\begin{remark} In the coordinate system $z\equiv x+\mu$, we have that $Z_1=(z^3-3\mu^2z)\partial/\partial z$ and $Z_2=z\partial/\partial z$. Hence, every Lie system of the form $Z=a_1(t)Z_1+a_2(t)Z_2$ can be written as a Bernoulli equation. As a consequence, the Abel equations satisfying the condition (\ref{CA}) can easily be integrated.
\end{remark}
Let us study the second type of quasi-Lie systems for $S(W_{\rm Abel},V_{\rm Abel})$,  namely those $X\in \mathcal{V}_{\rm Abel}$ which can be mapped onto
a Lie system related to a one-dimensional Vessiot--Guldberg Lie algebra through a curve  $g\in G_{W_{\rm Abel}}$. In other words, there exists a curve $g\in G_{W_{\rm Abel}}$ such that $g_\di X=\xi(t)(\sum_{\alpha=0}^3c_\alpha Y_\alpha)$ for certain constants $c_0,\ldots,c_3\in\mathbb{R}$ and a $t$-dependent function $\xi(t)$. That is, we have that
\begin{equation}\label{Transsys2}
\begin{aligned}
c_3\xi(t)&=f_3(t)\alpha^2(t),\\
c_2\xi(t)&=\alpha(t)(3f_3(t)\beta(t)+f_2(t)),\\
c_1\xi(t)&=3f_3(t)\beta^2(t)+2f_2(t)\beta(t)+f_1(t)-\frac{\dot\alpha(t)}{
\alpha(t)},\\
c_0\xi(t)&=\frac{1}{\alpha(t)}
\left(f_3(t)\beta^3(t)+f_2(t)\beta^2(t)+f_1(t)\beta(t)+f_0(t)-\dot\beta(t)\right
).\\
\end{aligned}
\end{equation}
Since $\alpha(t)\neq 0$, then $c_3\neq 0$. As a consequence, $\xi=f_3\,\alpha^2/c_3$ and
 we see that
\begin{equation}\label{sys4}
\begin{aligned}
\frac{c_2f_3\alpha}{c_3}&={3f_3\beta+f_2},\\
\frac{c_1f_3\alpha^2}{c_3}&=3f_3\beta^2+2f_2\beta+f_1-\frac{\dot
\alpha}{\alpha},\\
\frac{c_0f_3\alpha^2}{c_3}&=\frac{1}{\alpha}
\left(f_3\beta^3+f_2\beta^2+f_1\beta+f_0-\dot\beta\right
).
\end{aligned}
\end{equation}
Assuming that $c_2\neq 0$, we have that $3f_3\beta+f_2\neq 0$ and the first expression of the preceding system allows us to obtain \begin{equation}\label{l3}
\begin{gathered}
\frac{c_1c_3(3f_3\beta+f_2)^2}{c_2^2f_3}=3f_3\beta^2+2f_2\beta+f_1+\frac{f_2\dot f_3-f_3\dot f_2-3f_3^2\dot \beta}{f_3(f_2+3f_3\beta)},\\
\frac{c_3^2c_0(3f_3\beta+f_2)^3}{c_2^3f_3^2}=f_3\beta^3+f_2\beta^2+f_1\beta+f_0-\dot\beta.\\
\end{gathered}
\end{equation}
From the second equation, we can express the value of $\dot \beta$ in terms of $\beta$, $f_0,\ldots,f_3$ and the constants $c_0,\ldots,c_3$. Substituting the corresponding expression of $\dot \beta$ in the first one, we obtain that  $\beta(t)$  is a solution of a third-degree polynomial with respect to $\beta(t)$ and coefficients given by $f_0(t),\ldots,f_3(t)$ and $c_0,\ldots,c_3$. This shows that
\begin{equation}\label{formb}
\beta(t)=F(f_0(t),\ldots,f_3(t), c_0,\ldots, c_3),
\end{equation}
for a certain function $F$. This $\beta$ gives rise to a solution of system (\ref{l3}) if and only if it satisfies the second differential equation in (\ref{l3}). Since $\beta$ is of the form (\ref{formb}), this happens when the functions $f_0,\ldots,f_3$, and the constants $c_0,\ldots,c_3$, determining the transformed Abel equation obey a certain condition.
\begin{theorem}
An Abel equation of first-order and first-type  is a quasi-Lie system with respect to the scheme $S(W_{\rm Abel},V_{\rm Abel})$ and a one-dimensional Lie algebra $V=\langle c_3\partial/\partial x+c_2\partial/\partial x+c_1\partial/\partial x+c_0\rangle$, with $c_3\,c_2\neq 0$, if and only if system (\ref{l3}) admits a particular solution $\beta$. Given such a solution, we have that
\begin{align*}
\alpha=c_3\frac{3f_3\beta+f_2}{c_2f_3}.
\end{align*}
This gives the transformation $x(t)=\alpha(t)\bar x(t)+\beta$ mapping the first Abel equation onto the Lie system with Vessiot--Guldberg Lie algebra $V$.
\end{theorem}

\section{Integrability of some Abel equations}

\begin{theorem}
The knowledge of a particular solution of an Abel equation allows us to transform it into the canonical form
\begin{equation}\label{canAbel}
\frac{dx}{d\tau}=x^3+f_2(\tau)x^2.
\end{equation}
\end{theorem}
\begin{proof}
Using (\ref{NewCoff}), we see that if we construct a $t$-dependent change of variables $x(t)=\alpha(t)\bar x(t)+\beta(t)$ with $\beta(t)$ being a particular solution of the Abel equation (\ref{Abel}), then transformed Abel equation has, in view of (\ref{eqv4}), the $t$-dependent function $\bar{f}_0(t)=0$. From (\ref{NewCoff}) we see also that $\bar{ f}_1(t)$ becomes equal to zero for $\alpha(t)$ satisfying the easily integrable differential equation:
\begin{align*}
\frac{d\log \alpha}{dt}=3f_3(t)\beta^2(t)+2f_2(t)\beta(t)+f_1(t).
\end{align*}
Hence, the above $t$-dependent change of variables maps (\ref{Abel}) onto
\begin{align*}
\frac{d\bar x}{dt}=f_3(t)\alpha^2(t)\bar x^3+\alpha(t)(f_2(t)+3f_3(t)\beta(t))\bar x^2.
\end{align*}
Finally, the $t$-parametrisation
$
\tau=\int^tf_3(t')\alpha^2(t')dt'
$
transforms this equation into the canonical form (\ref{canAbel}).
\end{proof}

\begin{corollary} The Abel equation (\ref{canAbel}) is integrable for
$
9\dot f_2+2f_2^3=0.
$
\end{corollary}
\begin{proof} Under the assumed condition and in view of Theorem \ref{th1}, Abel equation (\ref{canAbel}) becomes a quasi-Lie system and it can be mapped onto a Lie system related to a bidimensional Vessiot--Guldberg Lie algebra. The change of variables $z\equiv x+\mu$ maps this Lie system onto a Bernoulli equation that can easily be solved.
\end{proof}

When the transformation property of our quasi-Lie scheme can be used to relate an
equation of the form (\ref{Abel}) with an integrable one, the initial equation can be solved.
It can be proved that we can use our quasi-Lie scheme for Abel equations to study also the integrability of Riccati and linear differential equations. The key point is mapping Eq.
(\ref{finalAbel})  onto an integrable Lie system
by means of a $t$-dependent change of variables induced by a curve $g$ in the group
$G_{W_{{\rm Abel}}}$ of $S(W_{{\rm Abel}},V_{{\rm Abel}})$. In this way, some results of
\cite{CarRam} are recovered.

\section{Quasi-Lie invariants}
Consider the Abel equation of first-order and first-kind related to the vector field
$X=\sum_{\alpha=0}^3A_\alpha(t)Y_{\alpha}$ with $Y_0,\ldots,Y_3$ given by (\ref{VectAbel}). We associate $X$ with a $t$-dependent function $F_{\rm Abel}(X)=\Phi_3^5/\Phi_5^3$, where
\begin{align*}
\Phi_3=\frac{dA_3}{dt}A_2-A_3\frac{dA_2}{dt}-3A_0A_3^2+A_1A_2A_3-\frac 29A_2^3,\\
 \Phi_5=-A_3\frac{d\Phi_3}{dt}-3\left(-\frac{dA_3}{dt}+\frac 13A_2^2-A_1A_3\right)\Phi_3.
\end{align*}
The map $F_{\rm Abel}:X\in \mathcal{V}_{\rm Abel}\mapsto F_{\rm Abel}(X)\in C^\infty(\mathbb{R})$ is such that $F_{\rm Abel}(X)=F_{\rm Abel}({g_\di X})$ for every $g\in\mathcal{G}_{W_{\rm Abel}}$ and $X\in \mathcal{V}_{\rm Abel}$, i.e.  $F_{\rm Abel}$ takes the same value on the elements of $\mathcal{V}_{\rm Abel}$ connected by an element of $\mathcal{G}_{W_{\rm Abel}}$ \cite[pg. 181]{Sw}. So, if two Abel equations associated with the $t$-dependent vector fields $X$ and $Y$ can be connected through an element of $\mathcal{G}_{W_{\rm Abel}}$, the function $F_{\rm Abel}$ must be the same for both. This can be employed to easily determine when two Abel equations are not connected.

In this section, we aim to use the theory of quasi-Lie schemes to study invariant functions like $F_{\rm Abel}$. As a result, we devise a method to easily obtain in a geometric way such invariants and to generalise such results to other systems of differential equations.

\begin{definition}\label{Def:QLI}
A {\it quasi-Lie invariant} relative to $\mathcal{G}_{W}$ is an $F\in {\rm Map}(\mathcal{V},C^\infty(\mathbb{R}))$ invariant under the action of the elements of $\mathcal{G}_{W}$, i.e.
\begin{align*}
F(g_\di X)=F(X),\qquad \forall g\in \mathcal{G}_{W},\qquad \forall X\in \mathcal{V}.
\end{align*}
\end{definition}
We see from Definition \ref{Def:QLI} that $F_{\rm Abel}$ is a quasi-Lie invariant relative to $S(W_{\rm Abel},V_{\rm Abel})$.
\begin{remark} To simplify the notation, we omit referring specifically to $\mathcal{G}_W$ when the group $G_W$ is understood from context.
\end{remark}

\begin{lemma} The space $\mathcal{I}(W,V)$ is a vector space relative to the standard sum of functions and their multiplication by real numbers. The space $(\mathcal{I}(W,V),\star)$  becomes an $\mathbb{R}$-algebra with unity relative to the product
\begin{equation}\label{star}
[(F\star G)(X)](t)\equiv [F(X)](t)[G(X)](t),\qquad F,G\in \mathcal{I}(W,V),\qquad \forall t\in\mathbb{R}.
\end{equation}
\end{lemma}

The following proposition, whose proof is straightforward, describes an operator that allows us to obtain new quasi-Lie invariants in terms of known ones.

\begin{proposition} The operator $D:\mathcal{I}(W,V)\rightarrow \mathcal{I}(W,V)$ of the form
\begin{align*}
[(DF)(X)](t)\equiv\frac{d}{dt} [F(X)(t)],\qquad \forall t\in\mathbb{R},
\end{align*}
is a derivation of the $\mathbb{R}$-algebra $(\mathcal{I}(W,V),\star)$.
\end{proposition}

Since every Lie system $X$ is described by a $t$-dependent vector field taking values in a Vessiot--Guldberg Lie algebra $V$, there always exists a $t$-dependent change of variables mapping $X$ onto another Lie system $\bar X$ taking values in $V$. As a consequence, we can prove the next proposition showing 
that Lie systems only admit trivial quasi-Lie invariants. 
  

\begin{proposition} If a quasi-Lie scheme $S(V,V)$ on a manifold $N$ admits a quasi-Lie invariant $F$, then $F(X)=F(\bar X)$ for every $X,\bar X\in \mathcal{V}$. 
\end{proposition}

A quasi-Lie invariant $F$ is said to be {\it constant} when $F(X)$ is the same for each $X\in \mathcal{V}$. %

\begin{example} Each Riccati equation is described by a $t$-dependent vector field $X$ taking values in the Vessiot--Guldberg Lie algebra $V_{\rm Ricc}$ given by (\ref{VGRicc}). We know that $S(V_{\rm Ricc},V_{\rm Ricc})$ is a quasi-Lie scheme. Every $X\in \mathcal{V}_{\rm Ricc}$ can be mapped onto any other Riccati equations through the action of a curve in $SL(2,\mathbb{R})$ on $\mathcal{V}_{\rm Ricc}$ (cf. \cite{Dissertationes}). Hence, as $F$ takes by definition the same value on those $X$ connected by a curve of $SL(2,\mathbb{R})$, then $F$ has the same value on the whole $\mathcal{V}_{\rm Ricc}$ and becomes constant. 
\end{example}

Subsequently, we focus on  the hereafter called quasi-Lie invariants of order $p$. We say that $F:X\in \mathcal{V}\mapsto F(X)\in C^\infty(\mathbb{R})$ is {\it a quasi-Lie invariant of order $p$ for $S(W,V)$} if there exists a function $\varphi_F:{\rm T}^pV\rightarrow \mathbb{R}$ such that
\begin{align*}
[F(X)](t)=\varphi_F\left({\bf t}_t^p X\right),\qquad \forall t\in\mathbb{R},\,\,\forall X\in \mathcal{V},
\end{align*}
where ${\bf t}_t^pX$ is the $p$-order tangent vector at $t$ of the curve $X\in \mathcal{V}$ and ${\rm T}^pV$ is the $p$-order tangent space to $V$, namely the manifold of all ${\bf t}_t^pX$ for a fixed $t$. 

To establish a coordinate system
on ${\rm T}^pN$, for $N$ being a manifold,  we need to use a connection $\nabla$. Nevertheless, as we aim to study ${\rm T}^pV$, for $V$ being a vector space, we may just assume $\nabla_{\partial/\partial y_\alpha}\partial/\partial y_\beta=0$,  where  $y_1,\ldots,y_r$ are the coordinate system with respect to the basis $\{X_1,\ldots,X_r\}$ of $V$, to obtain
\begin{align*}
{\bf t}^p_tX=\left(X_t,\frac{dX_t}{dt},\nabla_{\frac{dX_t}{dt}}\frac{dX_t}{dt},\ldots,\nabla^{p-1}_{\frac{dX_t}{dt}}\frac{dX_t}{dt}\right)=
\left(X_{t_0},\frac{dX_{t_0}}{dt},\frac{d^2X_{t_0}}{dt^2},\ldots, \frac {d^pX_{t_0}}{dt^p}\right),
\end{align*}
where $\frac {d^pX_{t_0}}{dt^p}\equiv\frac {d^pX_{t}}{dt^p}(t_0)$.
Note that if two curves $X_t$ and $Y_t$ have the same ${\bf t}_t^pX={\bf t}_t^pY$, then 
\begin{align*}
\left(X_t,\frac{dX_t}{d\bar {t}},\frac{d^2X_t}{d\bar{t}^2},\ldots, \frac {d^pX_t}{d\bar {t}^p}\right)=\left(X_t,\frac{dY_t}{d\bar t},\frac{d^2Y_t}{d\bar{t}^2},\ldots, \frac {d^pY_t}{d\bar{t}^p}\right)
\end{align*}
for every possible reparametrisation $\bar t=\bar t(t)$. Moreover, the coordinates of ${\bf t}^k_tX$ for a certain parameter $t$ depend linearly on the coordinates of ${\bf t}^k_tX$ for $\bar t$. Consequently, if ${\bf t}^k_tX=0$ for a certain $t$, it is so for every other reparametrisation $\bar t=\bar t(t)$.

If $\varphi_F$ exists, then it is necessarily unique. Moreover, our assumption involves that the value of $[F(X)](t_0)$ is determined exclusively by ${\rm t}^p_{t_0}X$. In other words, $[F(X)](t_0)$ depends on the values of $X$ on an interval of $t_0$. It is obvious, that there may exist other types of functionals $F$ that do not hold this property, e.g. $[\widetilde F(x)](t_0)\equiv \varphi_X({\bf t}^p_{t_1}X)$ for $t_1\neq t_0$.

We write $\mathcal{I}^{p}(W,V)$ for the space of quasi-Lie invariants of order $p$ for $S(W,V)$. Obviously, each $\mathcal{I}^p(W,V)$ and $\mathcal{I}(W,V)$ are $\mathbb{R}$-algebras relative to the product (\ref{star}).
As an instance, observe that the mapping $F_{\rm Abel}:X\in \mathcal{V}_{\rm Abel}\mapsto \Phi_3^5/\Phi_5^3 \in C^\infty(\mathbb{R})$ is a quasi-Lie invariant of order two for the quasi-Lie scheme $S(W_{\rm Abel},V_{\rm Abel})$. Let us study the properties of $p$-order quasi-Lie invariants and apply our results to the analysis of Abel equations.
 \begin{proposition} We have the following sequence of derivations
  \begin{align*}
  (\mathcal{I}^i(W,V),\star)\stackrel{D}{\longrightarrow}  (\mathcal{I}^{i+1}(W,V),\star),\qquad i\in \mathbb{N}.
    \end{align*}
 \end{proposition}

 \begin{theorem}\label{qlI} Let $S(W,V)$ be a quasi-Lie scheme. If  $F:\mathcal{V}\rightarrow C^\infty(\mathbb{R})$ is a quasi-Lie invariant of order $p$ for $S(W,V)$, then $\varphi_F$ is a common first-integral of the Lie algebra  $V_{\mathcal{J}^p}$ and $V_{\mathcal{T}^p}$ of fundamental vector fields corresponding to the Lie group actions:
\begin{align}
&\begin{array}{rccc}
\mathcal{J}^p:\!&G\times T^pV&\longrightarrow& T^pV\\
&(g;{\bf t}^p_{t_0}X=(X_0,\ldots,X_p))&\mapsto& {\bf t}^p_{t_0}(g_*X)=(g_*X_0,\ldots,g_*X_p),
\end{array}\,\label{action1}\\
&
\begin{array}{rccc}
\mathcal{T}^p:&T^p_{p-1}W\times T^pV&\longrightarrow& T^pV\\
&({\bf t}_{t_0}^pY;{\bf t}^p_{t_0}X=(X_0,\ldots,X_p))&\mapsto& {\bf t}^p_{t_0}(Y+X),
\end{array}\label{action2}
\end{align}
where $G$ is a connected Lie group with Lie algebra $T_eG\simeq W$, and $T^p_{p-1}W\simeq V$ is the Abelian Lie group spanned by those ${\bf t}^p_{t_0}Y\in T^pW$ such that ${\bf t}_{t_0}^{p-1}Y=0$. 
\end{theorem}

\begin{proof}

If $F$ is a quasi-Lie invariant, then $F(X)=F(g_\di X)$ for every $X\in \mathcal{V}$ and $g\in\mathcal{G}_W$. As $F$ is also a quasi-Lie invariant of order $p$, we get
\begin{align*}
\!\varphi_F\!\left({\bf t}_t^pX	\right)\!=\!\varphi_F\!\left({\bf t}_t^pg_\di X\right)\!=\varphi_F\!\left[{\bf t}^p_t\left(\frac{\partial g_t}{\partial t}\circ g_t^{-1}\!+\!g_{t*}X_t\right)\right]\!,
\end{align*}
for every $t\in\mathbb{R}$, $X\in \mathcal{V}$ and $g\in\mathcal{G}_W$. In particular, for a curve in $\mathcal{G}_W$ of the form $g:t\in \mathbb{R}\mapsto g^0\in G$, with $g^0$ a fixed element of $G$, the above expression becomes
\begin{align*}
\!\varphi_F\!\left({\bf t}_t^pX	\right)\!&=\varphi_F\!\left({\bf t}^p_tg^0_*X\right)=\!\varphi_F\!\left(g^0_{*}X_t,g^0_{*}\frac{dX_t}{dt},g^0_{*}\frac{d^2X_t}{dt^2},\ldots,g^0_{*}\frac{d^pX_t}{dt^p}\right)\!
=\!\varphi_F\left(\mathcal{J}^p\left(g^0,{\bf t}^p_tX\right)\right),
\end{align*}
for all $t\in\mathbb{R}$ and $X\in \mathcal{V}$. Hence, $\varphi_F$ is invariant along the orbits of $\mathcal{J}^p$, which amounts to saying that $\varphi_F$ is a common first integral for the elements of $V_{\mathcal{J}^p}$.

Let us now consider a curve $g^Y:t\in\mathbb{R}\mapsto e^{\int^t_{t_0} v_{t'}dt' }\in G$, where $e^{v}$ stands for the exponential of $v\in T_eG$, we write $v_t=f(t)v$ and ${\bf t}^{p-1}_{t_0}f=0$. Let $Y=(\partial {g}^Y/\partial t)\circ {g^Y}^{-1}$. Then,
\begin{align*}
\begin{gathered}
(g^Y_{\di} X)_t=Y_t+e^{\int^t_{t_0}v_{t'}dt'}_*X_t\Longrightarrow (g^Y_{\di} X)_{t_0}=X_{t_0},\\
\frac{d}{dt}(g^Y_{\di} X)_t=\frac{dY_t}{dt} -e^{\int^t_{t_0}v_{t'}dt'}_*\left(\left[Y_t,X_t\right]-\frac{dX_t}{dt}\right)\Longrightarrow \frac{d}{dt}(g^Y_{\di} X)\bigg|_{t_0}=\frac{dX_t}{dt}\bigg|_{t=t_0},\\
\frac{d^2}{dt^2}\bigg|_{t=t_0}(g^Y_{\di} X)_t=\frac{d^2Y_t}{dt^2} +e^{\int^t_{t_0}v_{t'}dt'}_*\left(\left[Y_t,\left[Y_t,X_t\right]\right]-\left[\frac{dY_t}{dt},X_t\right]\right)\\
-e^{\int^t_{t_0}v_{t'}dt'}_*\left(\left[Y_t,\frac{dX_t}{dt}\right]-\frac{d^2X_t}{dt^2}\right)\Longrightarrow \frac{d^2}{dt^2}\bigg|_{t=t_0}(g^Y_{\di} X)_t=\frac{d^2X_t}{dt^2}\bigg|_{t=t_0}.
\end{gathered}
\end{align*}
In view of the above and since ${\bf t}^{p-1}_{t_0}f=0$, the $d^{p}Y_t/dt^{p}$ is the first derivative of $Y$ that does not necessarily vanish. Therefore, we have
\begin{align*}
\frac{d^k}{dt^k}\bigg|_{t=t_0}(g^Y_{\di} X)_t=\frac{d^kX_t}{dt^k}\bigg|_{t=t_0},\quad k<p,\quad \frac{d^p}{dt^p}\bigg|_{t=t_0}(g^Y_{\di} X)_t=\frac{d^{p}Y_t}{dt^{p}}\bigg|_{t=t_0}+\frac{d^pX_t}{dt^p}\bigg|_{t=t_0}.
\end{align*}As a consequence,
\begin{align*}
\varphi_F\left({\bf t}^p_{t_0}X\right)=\varphi_F\left({\bf t}^p_{t_0}(X+Y)\right)=\varphi_F\circ \mathcal{T}^p\left({\bf t}^{p}_{t_0}Y,{\bf t}^{p}_{t_0}X\right),
\end{align*}
for all $Y\in W$, $X\in \mathcal{V}$. So, $\varphi_F$ is a common first integral of the vector fields of the Lie group action $\mathcal{T}^p$.
\end{proof}

\begin{remark} Observe that the above condition is purely geometrical, i.e. it does not depend neither on the coordinate system on $V$ nor on the variable $t$ on $\mathbb{R}$.
\end{remark}
\section{Quasi-Lie invariants of order two}
Let us provide a criterion to characterise the existence and to derive quasi-Lie invariants up to second-order. Next, we prove that Abel equations admit only four functionally invariant quasi-Lie invariants of order two containing, as a particular instance, the invariant $\Phi_3^5/\Phi_5^3$ known by Liouville. The description of higher-order quasi-Lie invariants can be approached in the same way, but we left their specific study for a future work.

A quasi-Lie invariant $F$ of order two is an element of ${\rm Map}(\mathcal{V},C^\infty(\mathbb{R}))$ for which there exists a function $\varphi_F:T^2V\rightarrow\mathbb{R}$ satisfying that
\begin{equation}\label{qlI2}
F(X)(t)=\varphi_F\left({\bf t}^2_tX\right),\qquad \forall t\in\mathbb{R},\qquad X\in \mathcal{V}.
\end{equation}
Let us determine the conditions ensuring that $F$ is a quasi-Lie invariant of order two. A continuation, we assume $V$ to be endowed with a connection $\nabla$ establishing a diffeomorphism $T^2V\simeq V^3$, which allows us to use local coordinates in higher-order tangent bundles. Proofs of following lemmas are not displayed as they are immediate after some calculations.

\begin{lemma} The map $\mathcal{T}_2:T^2_0W\times T^2V\rightarrow T^2V$, with $T^2_0W\!=\!\{{\rm t}^2_{t_0}Y\!\!\in\! T^2W|Y_{t_0}\!=\!0\}$ and
\begin{multline}\label{TL2}
\mathcal{T}_2\left({\bf t}^2_{t_0}Y=\left(0,\frac{dY_{t_0}}{d\bar t},\frac{d^2Y_{t_0}}{d\bar t^2}\right);{\bf t}_{t_0}^2X=\left(X_{t_0},\frac{dX_{t_0}}{d\bar t},\frac{d^2X_{t_0}}{d\bar t^2}\right)\right)\equiv\\ \left(X_{t_0},\frac{dY_{t_0}}{d\bar t}+\frac{dX_{t_0}}{d\bar t},-\frac{dt}{d\bar t}\left[\frac{dY_{t_0}}{d\bar t},X_{t_0}\right]+\frac{d^2X_{t_0}}{d\bar t^2}+\frac{d^2Y_{t_0}}{d\bar t^2}\right),
\end{multline}
is a Lie group action relative to the Abelian Lie group structure of $T^2_0W$.
\end{lemma}
\begin{lemma} Given a vector field $L\in W$, the mapping $\mathcal{T}_1^L:(s;{\bf t}^2_{t_0}X=(X,Y,Z))\in\mathbb{R}\times T^2V\mapsto \mathcal{T}_1^L(s;{\bf t}^2_{t_0}X)\in T^2V$ of the form
\begin{multline}\label{TL1}
\mathcal{T}_1^L(s;{\bf t}^2_{t_0}X)\!\equiv\! {\bf t}^2_{t_0}X\!+\!\left(\!sL,\!-s\frac{dt}{d\bar t}[L,X],s^2\left(\frac{dt}{d\bar t}\right)^2\!\![L,[L,X]]\!-[2s\frac{dt}{d\bar t}L,Y]\!-\!\frac{d^2t}{d\bar t^2}[sL,X]\!\right)
\end{multline}
is a Lie group action relative to $(\mathbb{R},+)$.
\end{lemma}


\begin{theorem}\label{The:QLI} 
The mapping $F:\mathcal{V}\rightarrow C^\infty(\mathbb{R})$ is a quasi-Lie invariant of order two for $S(W,V)$ if and only if $\varphi_F$ is a common first-integral of the Lie algebras of fundamental vector fields  $V_{\mathcal{J}^2}$ and $V_{\mathcal{T}^2}$ corresponding to the actions (\ref{action1}), (\ref{action2}) and the vector fields related to the flows $\mathcal{T}_1^L$
with $L\in W$ and $\mathcal{T}_2$.
\end{theorem}
\begin{proof}

In view of (\ref{qlI2}), determining the conditions ensuring that $F(X)=F(g_\di X)$ for every $X\in \mathcal{V}$ and $g\in\mathcal{G}_W$ amounts to proving that
$
\varphi_F\left({\bf t}^2_{t_0}g_\di X\right)=\varphi_F\left({\bf t}^2_{t_0}X\right),
$ for every $t\in\mathbb{R}$, $X\in \mathcal{V}$ and $g\in\mathcal{G}_W$. Let us analyse the above condition for four particular types of curves in $G_W$. This will give us some necessary conditions for $\varphi_F$ to give rise to a quasi-Lie invariant of order two.

Consider the curve $g^{sY}_t=e^{-\int^t_{t_0}svdt'}$ with $v\in T_eG_W$ and let $Y$ be the fundamental vector field associated with $v$.   In view of (\ref{Descom}), we have
\begin{align*}
(g^{sY}_{\di} X)_t&=sY+g^{sY}_{t*} X_t\Longrightarrow (g^Y_{\di} X)_{t_0}= sY+X_{t_0},\\
\frac{d}{d\bar t}(g^{sY}_{\di} X)_t&= -g^{sY}_{t*}\left(\frac{dt}{d\bar t}[sY,X_t]-\frac{dX_t}{d\bar t}\right)\Longrightarrow \frac{d}{dt}(g^{sY}_{\lambda \di} X)_{t_0}=-\frac{dt}{d\bar t}[ sY,X_{t_0}]+\frac{dX_{t_0}}{d\bar t},\\
\frac{d^2}{d\bar{t}^2}(g^{sY}_{\di} X)_t&= g^{sY}_{t*}\left(\left(\frac{dt}{d\bar t}\right)^2\left[ sY,\left[sY,X_t\right]\right]-2\frac{dt}{d\bar t}\left[ sY,\frac{dX_t}{d\bar t}\right]+\frac{d^2X_t}{d^2\bar t}-\frac{d^2t}{d\bar t^2}[sY,X_t]\right).
\end{align*}
Hence
$$
\frac{d^2}{dt^2}(g^{sY}_{\di} X)_{t_0}=\left(\frac{dt}{d\bar t}\right)^2[sY,[sY,X_{t_0}]]-2\frac{dt}{d\bar t}\left[sY,\frac{dX_{t_0}}{d\bar t}\right]+\frac{d^2X_{t_0}}{d\bar{t}^2}-\frac{d^2t}{d\bar t^2}[sY,X_{t_0}].
$$
The invariance of $F$ under such a curve leads to
\begin{align}\label{FirstOrder}
\varphi_F\left({\bf t}_{t_0}^2X\right)=\varphi_F\left({\bf t}_{t_0}^2g_{\di}X\right)=\varphi_F\left(\mathcal{T}_1^Y\left(s,{\bf t}_{t_0}^2X\right)\right),\quad
\end{align}
for all $X\in \mathcal{V}$. Hence, $\varphi_F$ is invariant along the orbits of $\mathcal{T}_1^Y$, which amounts to saying that $\varphi_F$ is a common first integral for the elements of $V^{\mathcal{T}_1^Y}$.

Let us analyse the invariance of $F$ under the flow of a $t$-dependent vector field $Y\in \mathcal{W}$ taking values in a one-dimensional quasi-Lie family and such that  $Y_{t_0}=0$. Let  $g^Y$ be the curve in $G_W$ giving rise to the flow of $Y$. In view of (\ref{Descom}), we obtain
\begin{align*}
\begin{gathered}
(g^Y_{\di} X)_t=Y_{t}+g^Y_{t*}X_t\Longrightarrow (g^Y_{\di} X)_{t_0}=X_{t_0},\\
\frac{d}{d\bar t}(g^Y_{\di} X)_t=\frac{dY_t}{d\bar t} -\frac{dt}{d\bar t}g^{Y}_{t*}\left[Y_t,X_t\right]+g^Y_{t*}\frac{dX_t}{d\bar t}\Longrightarrow \frac{d}{d\bar t}(g^Y_{\di} X)_{t_0}=\frac{dY_{t_0}}{d\bar t}+\frac{dX_{t_0}}{d\bar t},\\
\end{gathered}
\end{align*}
and
\begin{equation*}
\begin{gathered}
\frac{d^2}{d\bar{t}^2}(g^Y_{\di} X)_t= \frac{d^2Y_t}{d\bar t^2}-\frac{d^2t}{d\bar t^2}g^Y_{t*}[Y_t,X_t]+\left(\frac{dt}{d\bar t}\right)^2g^Y_{t*}[Y_t,[Y_t,X_t]]-\frac{dt}{d\bar t}g_{t*}^Y\left[\frac{dY_t}{d\bar t},X_t\right]\\-\frac{dt}{d\bar t}g_{t*}^Y2\left[Y_t,\frac{dX_t}{d\bar t}\right]+g^Y_{t*}\frac{d^2X_t}{d^2\bar t}
\Rightarrow \frac{d^2}{dt^2}(g^Y_{\di} X)_{t_0}=-\frac{dt}{d\bar t}\left[\frac{dY_{t_0}}{d\bar t},X_{t_0}\right]+\frac{d^2(X_{t_0}+Y_{t_0})}{d\bar{t}^2}.
\end{gathered}
\end{equation*}
As a consequence,
\begin{align}\label{SecondOrder}
\varphi_F\left({\bf t}^2_{t_0}X\right)=\varphi_F\circ \mathcal{T}_2\left({\bf t}^2_Y,{\bf t}^2_{t_0}X\right),
\end{align}
for all ${\bf t}^2_{t_0}Y\in T^2_0W$. So, $\varphi_F$ is invariant under the vector field associated with the flow $\mathcal{T}_2$.

Theorem \ref{qlI} entails $\varphi_F$ is constant along the orbits of $\mathcal{J}^2$ and $\mathcal{T}^2$, which amounts to saying that $\varphi_F$ is a common first integral for the fundamental vector fields of $V_{\mathcal{J}^2}$ and $V_{\mathcal{T}^2}$.

The above results show  that if $F$ is a quasi-Lie invariant of order two, then $\varphi_F$ must be a common first integral of the elements of $V_{\mathcal{J}^2}$, $V_{\mathcal{T}^2}$, $V^{\mathcal{T}_2}$ and the vector fields associated with the flows $\mathcal{T}^L$ with $L\in W$.

Let us prove the converse, namely, if $F$ is of the form (\ref{qlI2}) and $\varphi_F$ is invariant under the actions and vector fields mentioned above, then $F$ is a quasi-Lie invariant of order two.

Recall that every element of $\mathcal{G}_W$ can be written in an open interval around $t_0$ in the form
$$
g=\exp(-\lambda_1(t)v_1)\times\ldots\times\exp(-\lambda_r(t)v_r).
$$
 Using the Taylor expansion for the functions $\lambda_\alpha(t)$ around $t=t_0$, we see that
\begin{multline}\label{deconstar}
g=\exp(-[\lambda^3_1(t)+\lambda^{2}_1(t-t_0)^2/2+\lambda^{1}_1(t-t_0)]v_1)\times\ldots\times\\ \exp(-[\lambda^3_r(t)+\lambda^{2}_r(t-t_0)^2/2+\lambda^{1}_r(t-t_0)]v_r)g_0,
\end{multline}
where $d^2\lambda_\alpha^3/dt^2(t_0)=d\lambda^3_\alpha(t)/dt(t_0)=\lambda^3_\alpha(t_0)=0$ for $\alpha=1,\ldots,r$ and $\lambda_i^j\in\mathbb{R}$ with $i=1,\ldots,r$ and $j=1,2$.
 Taking into account that $\varphi_F$ is invariant under the action $\mathcal{J}^2$, we get that  
\begin{align*}
\varphi_F\left({\bf t}^2_{t}g_{0\di}X\right)=\varphi_F\left(\mathcal{J}^2\left(g_0,{\bf t}^2_tX\right)\right)=\varphi_F\left({\bf t}^2_{t}X\right),\qquad \forall t\in\mathbb{R},
\end{align*}
for every $g_0\in G$ and $X\in \mathcal{V}$. Hence, $F(X)=F(g_{0\di}X)$ for all $g_0\in G$ and $X\in \mathcal{V}$.
Since $\varphi_F$ is invariant under the action of every $\mathcal{T}^Y_1$ and in view of (\ref{FirstOrder}), we have
\begin{multline}\label{con1}
F\left(e^{\lambda_\alpha^1(t_0-t)v_\alpha}_{\di}X\right)(t_0)= \varphi_F\left({\bf t}_{t_0}^2e^{\lambda_\alpha^1(t_0-t)v_\alpha}_{\di}X\right)=\\\varphi_F\circ \mathcal{T}^{Y_\alpha}_1\left(\lambda_\alpha^1,({\bf t}^2_{t_0}X)\right)=\varphi_F({\bf t}^2_{t_0}X)=F(X)(t_0).
\end{multline}
As $\varphi_F$ is invariant under the action of $\mathcal{T}_2$ and in view of (\ref{SecondOrder}), then
\begin{align}
F\left(e^{-\lambda_\alpha^2(t_0-t)^2v_\alpha/2}_{\di}X\right)(t_0)=\varphi_F\circ \mathcal{T}_2\left(0,\lambda_\alpha^2Y_\alpha,0),{\bf t}_{t_0}^2X\right)=\varphi_F\left({\bf t}_{t_0}^2X\right)=F(X)(t_0).
\end{align}
Similarly, as $\varphi_F$ is invariant under $\mathcal{T}_2$, we obtain that
\begin{align}\label{con3}
F\left(e^{-\lambda_\alpha^3(t)v_\alpha}_{\di}X\right)(t_0)=\varphi_F\circ \mathcal{T}_2\left((0,0,d^3\lambda_\alpha^3/dt^3(t_0)Y_\alpha),{\bf t}_{t_0}^2X\right)=\varphi_F\left({\bf t}_{t_0}^2X\right)=F(X)(t_0).
\end{align}
Using (\ref{con1})-(\ref{con3}) and decomposition (\ref{deconstar}), we get that $F(X)(t_0)=F(g_\di X)(t_0)$. This holds for any point $t_0$. So, $F(X)=F(g_\di X)$ and $F$ is a quasi-Lie invariant of order two.

\end{proof}

\section{Quasi-Lie invariants of first-order for Abel equations}
Let us prove that Abel equations admit, apart from the quasi-Lie invariant developed by Liouville, other three functionally independent ones.

The quasi-Lie invariants for Abel equations of order two must satisfy that their corresponding $\varphi_F$ must be a first-integral of:
\begin{enumerate}
 \item The fundamental vector fields of the actions $\mathcal{T}^2$ and $\mathcal{J}^2$ given by (\ref{action1})--(\ref{action2}) for $p=2$.
 \item The vector fields corresponding to the flows  $\mathcal{T}^L_1$ with $L\in W_{\rm Abel}$ given by (\ref{TL1}).

 \item The fundamental vector fields associated with $\mathcal{T}_2$ given by (\ref{TL2}).
\end{enumerate}

First, the fundamental vector fields of the Lie group action $\mathcal{J}^2$ are spanned by
\begin{align}\label{Ve1}
Y^{\mathcal{J}^2}_1=\sum_{\alpha=0}^2\left(\lambda^{\alpha)}_1\frac{\partial}{\partial \lambda^{\alpha)}_0}+2\lambda^{\alpha)}_2\frac{\partial}{\partial \lambda^{\alpha)}_1}+3\lambda^{\alpha)}_3\frac{\partial}{\partial \lambda^{\alpha)}_2}\right),\\ Y^{\mathcal{J}^2}_2=\sum_{\alpha=0}^2\left(-\lambda^{\alpha)}_0\frac{\partial}{\partial \lambda^{\alpha)}_0}+\lambda^{\alpha)}_2\frac{\partial}{\partial \lambda^{\alpha)}_2}+2\lambda^{\alpha)}_3\frac{\partial}{\partial \lambda^{\alpha)}_3}\right).
\end{align}
Since $[Y^{\mathcal{J}^2}_1,Y^{\mathcal{J}^2}_2]=-Y^{\mathcal{J}^2}_1$, then $\langle Y^{\mathcal{J}^2}_1,Y^{\mathcal{J}^2}_2\rangle $ is a non-Abelian two-dimensional Lie algebra.

Subsequently, the fundamental vector fields of the action $\mathcal{T}^2$ are linear combinations of 
\begin{align}\label{Ve2}
Y_1^{\mathcal{T}^2}=\frac{\partial}{\partial \lambda^{2)}_0},\qquad Y_2^{\mathcal{T}^2}=\frac{\partial}{\partial \lambda^{2)}_1},
\end{align}
which generate a two-dimensional Abelian Lie algebra.

The commutators between the elements of the bases of the two above Lie algebras read:
\begin{align*}
[Y_1^{\mathcal{J}^2},Y_1^{\mathcal{T}^2}]=0,\,\,[Y_1^{\mathcal{J}^2},\,\,Y_2^{\mathcal{T}^2}]=-Y_1^{\mathcal{T}^2},\,\,[Y_2^{\mathcal{J}^2},Y_1^{\mathcal{T}^2}]=Y_1^{\mathcal{T}^2},\,\,[Y_2^{\mathcal{J}^2},Y_2^{\mathcal{T}^2}]=0.
\end{align*}
So, they span a four-dimensional Lie algebra.

Every fundamental vector field of $\mathcal{T}_2(0,L,0),\cdot)$, with $L\in W_{\rm Abel}$, is a linear combination of 
\begin{align}\label{Ve3}
\Theta_2^{X_1}\!=\!\frac{\partial}{\partial \lambda^{1)}_0}\!-\!\lambda^{0)}_1\frac{\partial}{\partial \lambda^{2)}_0}\!-\!2\lambda^{0)}_2\frac{\partial}{\partial \lambda^{2)}_1}\!-\!3\lambda^{0)}_3\frac{\partial}{\partial\lambda^{2)}_2},
\,\, \Theta_2^{X_2}\!=\!\frac{\partial}{\partial \lambda^{1)}_1}\!+\!\lambda^{0)}_0\frac{\partial}{\partial \lambda^{2)}_0}\!-\!\lambda^{0)}_2\frac{\partial}{\partial \lambda^{2)}_2}\!-\!2\lambda^{0)}_3\frac{\partial}{\partial\lambda^{2)}_3}.
\end{align}
These vector fields hold $[\Theta_2^{X_1},\Theta_2^{X_2}]=0$. The non-zero commutation relations between these vector fields and all the previous ones read
$\left[\Theta_2^{X_1},Y_2^{\mathcal{J}_2}\right]=-\Theta_2^{X_1}
\left[\Theta_2^{X_2},Y_1^{\mathcal{J}_2}\right]=\Theta_2^{X_1}.$
  The flows $\mathcal{T}_1^{X_1}$ and $\mathcal{T}_1^{X_2}$ are related to the vector fields
\begin{align}\label{Ve4}
\Theta_1^{X_1}=\frac{\partial}{\partial \lambda^{0)}_0}-\lambda^{0)}_1\frac{\partial}{\partial \lambda^{1)}_0}-2\lambda^{0)}_2\frac{\partial}{\partial \lambda^{1)}_1}-3\lambda^{0)}_3\frac{\partial}{\partial\lambda^{1)}_2}-2\lambda^{1)}_1\frac{\partial}{\partial \lambda^{2)}_0}-4\lambda^{1)}_2\frac{\partial}{\partial \lambda^{2)}_1}-6\lambda^{1)}_3\frac{\partial}{\partial\lambda^{2)}_2},\nonumber\\
\Theta_1^{X_2}=\frac{\partial}{\partial \lambda^{0)}_1}+\lambda^{0)}_0\frac{\partial}{\partial \lambda^{1)}_0}-\lambda^{0)}_2\frac{\partial}{\partial \lambda^{1)}_2}-2\lambda^{0)}_3\frac{\partial}{\partial\lambda^{1)}_3}+2\lambda^{1)}_0\frac{\partial}{\partial \lambda^{2)}_0}-2\lambda^{1)}_2\frac{\partial}{\partial \lambda^{2)}_2}-4\lambda^{1)}_3\frac{\partial}{\partial\lambda^{2)}_3}.
\end{align}
The remaining vector fields associated with the flows $\mathcal{T}^L_1$ and $\mathcal{T}^L_2$ are linear combinations of the above ones. So, it is enough to find their common first-integrals to obtain the common first-integrals for the above ones. Meanwhile, we have that $[\Theta_1^{X_1},\Theta_1^{X_2}]=2\Theta_2^{X_1}$ and
the non-zero commutation relations between an element of $\{\Theta_1^{X_1},\Theta_1^{X_2}\}$ and an element of  $\{Y_1^{\mathcal{J}_2},Y_2^{\mathcal{J}_2},Y_1^{\mathcal{T}_2},Y_2^{\mathcal{T}_2}\}$ read
$$
\begin{array}{llllllllll}
\left[\Theta_1^{X_1},Y_2^{\mathcal{J}_2}\right]=-\Theta_1^{X_1},\,\, &\left[\Theta_1^{X_1},\Theta_2^{X_2}\right]=3Y_1^{\mathcal{T}_2}, \,\, &\left[\Theta_1^{X_2},Y_1^{\mathcal{J}_2}\right]=\Theta_1^{X_1},\,\, &\left[\Theta_1^{X_2},\Theta_2^{X_1}\right]=-3Y_1^{\mathcal{T}_2}.
\end{array}
$$

The above vector fields (\ref{Ve1})--(\ref{Ve4}) span an integrable distribution of rank eight on a manifold of dimension twelve. So, there must be at least four local integrals which give rise to quasi-Lie invariants.  One of them is the one known due to Liouville, i.e. $F=\Phi_3^5/\Phi_5^3$ for 
$$\begin{gathered}
\Phi_3=\lambda_3^{1)}\lambda_2^{0)}-\lambda_3^{0)}\lambda_2^{1)}-3\lambda_0^{0)}(\lambda_3^{0)})^2+\lambda^{0)}_1\lambda^{0)}_2\lambda^{0)}_3-2(\lambda^{0)}_2)^3/9,\\\Phi_5=-\lambda^{0)}_3 {\rm D}\Phi_3-3\left[-\lambda_3^{1)}+(\lambda^{0)}_2) ^2/3-\lambda_1^{0)}\lambda_3^{0)}\right]\Phi_3,\\
{\rm D}\Phi_3=\lambda_3^{2)}\lambda^{0)}_2-\lambda_3^{0)}\lambda^{2)}_2-3\lambda_0^{1)}(\lambda_3^{0)})^2-6\lambda_0^{ 0)}\lambda_3^{ 0)}\lambda_3^{1)}+\qquad\qquad\qquad\qquad\qquad\qquad\qquad\qquad\\\qquad\qquad\qquad\qquad\qquad\qquad\qquad\qquad+\lambda_1^{ 1)}\lambda^{0)}_2\lambda^{0)}_3+\lambda_1^{ 0)}\lambda^{1)}_2\lambda^{0)}_3+\lambda_1^{ 0)}\lambda^{0)}_2\lambda^{1)}_3-2(\lambda_2^{0)})^2\lambda_2^{1)}/3.
\end{gathered}
$$
Observe that $F$ does not depend on $\lambda_0^{2)}$, $\lambda_3^{2)}$ and $\lambda_1^{2)}$. So, $F$ is a common first-integral of $Y^\mathcal{T}_1$ and $Y^\mathcal{T}_2$.  It is easy to verify by means of a program of symbolic calculation that $F$ is a first-integral of all the remaining vector fields (\ref{Ve1})--(\ref{Ve4}).

Our geometric approach shows that there are more quasi-Lie integrals of order two apart from Liouville's one.
These invariants are also geometric: any other system of differential equations with an isomorphic algebraic quasi-Lie scheme have the same quasi-Lie invariants.

\section{Quasi-Lie invariants of order lower than two and Abel equations}

A quasi-Lie invariant of order zero (resp. one) can be considered  as a quasi-Lie invariant of order two whose $\varphi_X:T^2V\rightarrow \mathbb{R}$ depends only on $V$ (resp. $TV$). Hence, these quasi-Lie invariants can be characterised by restricting Theorem \ref{The:QLI} to these types of $\varphi_F$ as follows.  
\begin{corollary}\label{the:qlI0} A function $F:\mathcal{V}\rightarrow C^\infty(\mathbb{R})$ is a quasi-Lie invariant of order zero for $S(W,V)$ if and only if $\varphi_F$ is a common first-integral of the Lie algebras of fundamental vector fields  of the Lie group actions $\mathcal{J}^0$ and $\mathcal{T}^0$.
\end{corollary}

\begin{corollary}\label{the:qlI1} The mapping $F:\mathcal{V}\rightarrow C^\infty(\mathbb{R})$ is a quasi-Lie invariant of order one for $S(W,V)$ if and only if  $\varphi_F$ is a common first-integral of the Lie algebras of fundamental vector fields  of the Lie group actions ${\mathcal{J}^1}$, ${\mathcal{T}^1}=\mathcal{T}_1:T^1_0W\times TV\longrightarrow TV$, with
\begin{equation*}
\mathcal{T}_1\left({\bf t}_{t_0}Y=\left(0,\frac{dY_{t_0}}{d\bar t}\right);{\bf t}_{t_0}X=\left(X_{t_0},\frac{dX_{t_0}}{d\bar t}\right)\right)\equiv \left(X_{t_0},\frac{dX_{t_0}}{d\bar t}+\frac{dY_{t_0}}{d\bar t}\right),
\end{equation*}
with $T_0W=\{{\rm t}_{t_0}Y\in TW\,|Y(t_0)=0\}$, and the vector fields corresponding to the flows
\begin{align*}
\begin{array}{rccc}
\mathcal{T}^Z:&\mathbb{R}\times TV&\longrightarrow& TV\\
&(s;(X,Y))&\mapsto& (sZ+X,-\frac{dt}{d\bar t}[sZ,X]+Y),
\end{array}\qquad Z\in W.
\end{align*}
\end{corollary}

Let us apply the above theory to show that Abel equations do not admit smooth non-constant quasi-Lie invariants of orders zero and one. From the calculations of the above section, we get
$
{\rm t}^1_{t_0}X=\left(\left(\sum_{\alpha=0}^3\lambda^0_\alpha x^\alpha\right)\partial/\partial x, \left(\sum_{\alpha=0}^3\lambda^1_\alpha x^\alpha\right)\partial/\partial x\right)
$
and
\begin{equation*}
\mathcal{J}^1((b,a),{\rm t}^1_{t_0}X)=\Phi_{(b,a)*}{\rm t}^1_{t_0}X=\left(\sum_{\alpha=0}^3\lambda^0_\alpha (x-b)^\alpha a^{1-\alpha}\frac{\partial}{\partial x},\sum_{\alpha=0}^3\lambda^1_\alpha (x-b)^\alpha a^{1-\alpha}\frac{\partial}{\partial x}\right).
\end{equation*}
Hence, the Lie algebra of fundamental vector fields of $\mathcal{J}^1$ is spanned by
\begin{align*}
Y_1^{\mathcal{J}^1}=\sum_{\alpha=0}^1\left(-\lambda^\alpha_0\frac{\partial}{\partial \lambda^\alpha_0}+\lambda^\alpha_2\frac{\partial}{\partial \lambda^\alpha_2}+2\lambda^\alpha_3\frac{\partial}{\partial \lambda^\alpha_3}\right),\quad
Y_2^{\mathcal{J}^1}=\sum_{\alpha=0}^1\left(\lambda^\alpha_1\frac{\partial}{\partial \lambda^\alpha_0}+2\lambda^\alpha_2\frac{\partial}{\partial \lambda^\alpha_1}+3\lambda^\alpha_3\frac{\partial}{\partial \lambda^\alpha_2}\right).
\end{align*}
The Lie algebra of fundamental vector fields for $\mathcal{T}^1$ admits a basis 
$Y_1^{\mathcal{T}^1}={\partial}/{\partial \lambda^1_0}$, $Y_2^{\mathcal{T}^1}={\partial}/{\partial \lambda_1^1}$.
From ${\mathcal{T}}^Z:\mathbb{R}\times TV\rightarrow TV$, we obtain the vector field
\begin{equation*}
\mathcal{Z}=\frac{d}{ds}\bigg|_{s=0}\mathcal{T}^Z(s,(X,Y))=(Z,-\frac{dt}{d\bar t}[Z,X])\in { T}_{(X,Y)}TV.
\end{equation*}
Taking $Z=\mu_0\partial/\partial x+\mu_1x\partial/\partial x$, we have $\mathcal{Z}=-\mu_0\mathcal{Z}_1-\mu_1\mathcal{Z}_2$ with
\begin{align*}
\mathcal{Z}_1=-\frac{\partial}{\partial \lambda_0^0}+\lambda^0_1\frac{\partial}{\partial \lambda_0^1}+2\lambda^0_2\frac{\partial}{\partial \lambda_1^1}+3\lambda^0_3\frac{\partial}{\partial \lambda^1_2},\qquad
\mathcal{Z}_2=-\frac{\partial}{\partial \lambda^0_1}-\lambda^0_0\frac{\partial}{\partial \lambda_0^1}+\lambda^0_2\frac{\partial}{\partial \lambda^1_2}+2\lambda^0_3\frac{\partial}{\partial \lambda_3^1}.
\end{align*}
It can be shown that $[\mathcal{Z}_1,Y^2_\mathcal{J}]$ and $[\mathcal{Z}_2,Y^1_\mathcal{J}]$ and previous vector fields span a distribution of order eight on an open dense subset of $V^2\simeq \mathbb{R}^8$. So, their only common first-integrals are constant functions and the quasi-Lie invariants of order one are the trivial constants on $TV$.

\section{Conclusions and Outlook}

We have defined and studied the representation associated with a quasi-Lie scheme, algebraic morphisms of quasi-Lie schemes, structure preserving groups and Lie almost symmetries for quasi-Lie families. This retrieved the hierarchies of Lie systems devised in \cite{BM10}. We have also defined, studied and applied quasi-Lie invariants. We provided a geometric characterisation of some of these quasi-Lie invariants.

As an application, we proved that Abel equations can be characterised as elements of a quasi-Lie scheme. The transformations to study Abel equations have appeared as a consequence of the scheme. Using morphism of schemes, we generalised properties of Abel equations to other differential equations on the plane.  We found that the Liouville invariant is a quasi-Lie invariant for Abel equations and we prove the existence of some new ones.

In the future, we aim to continue studying quasi-Lie schemes and to apply our results to
Abel equations, Schr\"odinger equations and to study other types of quasi-Lie invariants.

\section*{Acknowledgements}
Partial financial support from MINECO (Spain) grant number MTM2012-33575. The research of J. de Lucas was
supported by the National Science Centre (POLAND)
under the grant HARMONIA No. 2012/04/M/ST1/00523.
J. de Lucas acknowledges a grant FMI40/10 from the Direcci\'on General de Arag\'on
to perform a research stay in the University of Zaragoza.

\end{document}